\documentstyle[12pt]{article}

\begin{document}
\medskip
\rightline {math.QA/9805032 \hskip 6mm}
\vskip 50 pt

\centerline{\bf REPRESENTATIONS OF THE QUANTUM
ALGEBRA $U_q({\rm u}_{n,1})$}
\vskip 15 pt

\centerline {\sc V. A. Groza, N. Z. Iorgov, and A. U. Klimyk}
\vskip 5 pt

\centerline {\it Institute for Theoretical Physics, Kiev 252143,
Ukraine}
\vskip 30 pt

\begin{abstract}
Infinite dimensional representations of the real form $U_q({\rm u}_{n,1})$
of the Drinfeld--Jimbo algebra $U_q({\rm gl}_{n+1})$ are defined. The
principal series of representations of $U_q({\rm u}_{n,1})$ is studied.
Intertwining operators for pairs of the principal series representations
are calculated in an explicit form. The structure of reducible
representations of the principal series is determined. Irreducible
representations of $U_q({\rm u}_{n,1})$, obtained from irreducible and
reducible principal series representations, are classified. All
$*$-representations in this set of irreducible representations
are separated. Unlike the classical case, the algebra
$U_q({\rm u}_{n,1})$ has finite dimensional irreducible $*$-representations.
\end{abstract}
\bigskip

\bigskip


\newpage

\noindent{\bf 1. Introduction}
\bigskip

\noindent
Under quantum algebras we mean quantized universal enveloping algebras
defined by Drinfeld [1] and Jimbo [2]. We denote them by
$U_q({g})$, where $q$ is a deformation parameter and ${ g}$ is
a corresponding semisimple Lie algebra.
Finite dimensional irreducible representations of $U_q({g})$
were classified by Rosso [3] and Lusztig [4] (see also [5] and [6]).
It was shown that every
irreducible finite dimensional representation of a simple Lie algebra
${g}$ can be deformed to be an irreducible representation of the
corresponding quantum algebra $U_q({g})$.
Moreover, if $q$ is not a root of unity, then every finite dimensional
irreducible representation of
$U_q({g})$ is essentially obtained in this way (after possibly
tensoring by a one-dimensional representation).

Along with finite dimensional representations of quantum algebras
their infinite dimensional representations are also important.
For example, infinite dimensional representations of the algebra
$U_q({\rm su}_{1,1})$ are closely related to the $q$-oscillator algebra.
They were used in the theory of $q$-special fucntions.
Our aim in this paper is to construct infinite dimensional representations
of the real form $U_q({\rm u}_{n,1})$ of the quantum algebra
$U_q({\rm gl}_{n+1})$.

Simplest infinite dimensional irreducible representations of quantum
algebras are determined by Verma modules. The theory of Verma modules of
$U_q({ g})$ when $q$ is not a root of unity is similar to that for the
corresponding universal enveloping algebras $U({g})$ (the details of this
construction see in [4] and [6]) and hence we do not concern
Verma modules in this paper.

It is possible to give different definitions of infinite dimensional
irreducible representations of $U_q({\rm u}_{n,1})$. We use the definition
similar to that for the classical case.
As in the classical case, we demand that the restriction of a representation
$T$ of $U_q({\rm u}_{n,1})$ to a maximal compact
subalgebra $U_q({\rm u}_n \oplus {\rm u}_1)$
decomposes into a direct sum of finite dimensional irreducible
representations of this subalgebra.

 In the case of infinite dimensional representations of real
semisimple Lie algebras, there is a one-to-one correspondence between
irreducible infinite dimensional representations of a connected simply
connected real semisimple Lie group $G$ and those of its Lie algebra $g$
(we consider only those representations of $G$ and $g$ which correspond to
Harish-Chandra $g$-modules). Let $P=MAN$ be a minimal parabolic subgroup in
$G$ [7], and let $\omega $ be an irreducible finite dimensional
representation of $P$. We induce by $\omega$ the representation $T_{\omega }$
of $G$. The representations $T_{\omega}$ constitute the principal
nonunitary series of $G$. Using the method described in [8] we can
construct the corresponding representations $T_{\omega}$ of $g$
induced by representations $\omega$ of the subalgebra $p=m+a+n$ (the
Lie algebra of $P$). By Harish-Chandra's theorem,
{\it every irreducible representation of $g$ (or of $G$) is equivalent to
some irreducible representation $T_{\omega}$ or to a subrepresentation
in a quotient representation of some reducible representation $T_{\omega}$.}
By means of this theorem a classification of all irreducible
representations of $g$ and of $G$ has been obtained.

    An essential item in construction of these representations is the
Iwasawa decomposition of $G$ (and of $g$) and parabolic subgroups.
This method cannot be extended to the quantum case since we have no
Iwasawa decomposition and parabolic subalgebras for real forms of
quantum algebras.

In order to construct infinite dimensional representations of
$U_q({\rm u}_{n,1})$ we write down explicit formulas for representation
operators corresponding to the generating elements of $U_q({\rm u}_{n,1})$.
Then we prove that they indeed determine representations of
$U_q({\rm u}_{n,1})$. In this way we obtain for $U_q({\rm u}_{n,1})$
a $q$-analogue of the principal nonunitary series. Then we succeed as
in the classical case, that is, we separate in this set of representations
irreducible ones and find all possible irreducible components of
reducible representations. Thus, we receive a set $\Omega$ of irreducible
representations of $U_q({\rm u}_{n,1})$. It will be shown in the
separate paper that the set $\Omega$ exhausts all irreducible
representations of $U_q({\rm u}_{n,1})$. By the standard method, irreducible
$*$-representations are separated in the set $\Omega$. Note that the
results of this paper were announced (without proofs) in the preprint [9].

Since in the set $\Omega$ there are irreducible finite dimensional
representations of $U_q({\rm u}_{n,1})$, then we obtain a new proof of
the Gel'fand--Tsetlin formulas for finite dimensional representations of the
quantum algebra $U_q({\rm gl}_{n+1})$. Remark that these formulas for
$U_q({\rm gl}_{n+1})$ were stated without proof by Jimbo [10]. The proof
was given in [11]. However, this proof is rather very complicated. Our
proof is much more simple and natural.
\medskip

\noindent
{\bf 2. The quantum algebras $U_q({\rm u}_n({\epsilon}))$}
\medskip

The quantum algebra $U_q({\rm gl}_n)$ is a $q$-deformation of the universal
enveloping algebra $U({\rm gl}(n,{\bf C}))$ and is
generated by the elements
$$k_i,\ \ k_i^{-1},\ \ e_j,\ \ f_j,\ \ \ \
i=1,2,\cdots ,n, \ \ j=1,2,\cdots ,n-1,$$
satisfying the relations
$$k_ik_i^{-1}=k_i^{-1}k_i=1,\ \ \ k_ik_j=k_jk_i, \eqno (1)$$
$$k_ie_jk_i^{-1}=q^{\delta _{ij}-\delta _{i,j+1}}e_j,\ \ \
k_if_jk_i^{-1}=q^{-\delta _{ij}+\delta _{i,j+1}}f_j,
\eqno(2)$$
$$[e_i,f_j]=\delta _{ij}{k_ik_{i+1}^{-1}-k^{-1}_ik_{i+1}\over q-q^{-1}} ,
\eqno(3)$$
$$[e_i,e_j]=[f_i,f_j]=0,\ \ \ \vert i-j\vert >1, \eqno(4)$$
$$e_i^2e_{i\pm 1}-(q+q^{-1})e_ie_{i\pm 1}e_i+e_{i\pm 1}e^2_i=0,
\eqno (5)$$
$$f_i^2f_{i\pm 1}-(q+q^{-1})f_if_{i\pm 1}f_i+f_{i\pm 1}f^2_i=0.
\eqno(6)$$
The last two relations are called the $q$-Serre relations. A structure of a
Hopf algebra is introduced into $U_q({\rm sl}_n)$ (see, for example [5]).
We do not need this Hopf algebra structure in this paper and shall use
only the algebraic structure of $U_q({\rm gl}_n)$.

We suppose throughout the paper that $q$ is a positive real number and
equip $U_q({\rm gl}_n)$ with $*$-structures. Different $*$-structures
define real forms of $U_q({\rm gl}_n)$. Let $\epsilon :=
(\epsilon _0,\epsilon _1,\cdots ,\epsilon _{n-1})$, $\epsilon _i=\pm 1$.
The $*$-structure of $U_q({\rm gl}_n)$ associated with $\epsilon$ is
determined by
$$(k_i^{\pm 1})^*=k_i^{\pm 1},\ \ \ e_i^*=\epsilon _{i-1}\epsilon _if_i,
\ \ \ f_i^*=\epsilon _{i-1}\epsilon _i e_i
$$
and is denoted by $U_q({\rm u}_n(\epsilon ))$. In particular,
the $*$-structure determined by
$$(k_i^{\pm 1})^*=k_i^{\pm 1},\ \ \ e^*_p=-f_p,\ \ \ f^*_p=-e_p, \eqno(7)$$
$$e^*_i=f_i,\ \ \ f^*_i=e_i,\ \ \ i\ne p, \eqno (8)$$
defines the real form denoted by
$U_q({\rm u}_{p,n-p})$. The real form
$U_q({\rm u}_n(\epsilon ))$ with $\epsilon =(1,1,\cdots$, $1)$ is called
the compact real form of $U_q({\rm gl}_n)$ and is denoted by
$U_q({\rm u}_n)$.

Note that the quantum algebras $U_q({\rm u}_n(\epsilon ))$
are of a considerable interest since they are closely related to the
quantum hyperboloids. The quantum hyperboloid $M^{2n-1}_q$ is
the associative algebra generated by the elements $z_i$, ${\hat z}_i$,
$i=0,1,\cdots ,n-1$, and $\kappa$ with the defining relations
$$
z_iz_j=q^{-1}z_jz_i,\ \ \ {\hat z}_i{\hat z}_j=q{\hat z}_j{\hat z}_i, \ \ \
i<j, $$
$$
z_i{\hat z}_j=q^{-1}{\hat z}_jz_i,\ \  i\ne j,\ \ \
{z}_i{\hat z}_i-{\hat z}_iz_i =(q^2-1)\Biggl( \sum _{k\ge i}
z_k{\hat z}_k+\kappa \Biggr) , $$
$$
z_i\kappa =q^{-2}\kappa z_i,\ \ \ {\hat z}_i\kappa =q^2\kappa {\hat z}_i
$$
and equipped with the $*$-structure determined by the formulas
$$
z_i^*=\epsilon _i \epsilon _{n-1}{\hat z}_i,\ \ \ \kappa ^*=\kappa .
$$
The algebra $U_q({\rm u}_n(\epsilon ))$ acts on $M_q^{2n-1}$ as
$$
e_i:\ z_j\to \delta _{ij}q^{1/2}z_{j-1},\ \ \
{\hat z}_j\to -\delta _{i-1,j}q^{3/2}{\hat z}_{j+1},\ \ \ \kappa \to 0, $$
$$
f_i:\ z_j\to \delta _{i-1,j}q^{-1/2}z_{j+1},\ \ \
{\hat z}_j\to -\delta _{ij}q^{-3/2}{\hat z}_{j-1},\ \ \ \kappa \to 0, $$
$$
k_i:\ z_{i-1}\to q^{1/2}z_{i-1},\ \ \
{\hat z}_{i-1}\to q^{-1/2}{\hat z}_{i-1},\ \ \ \kappa \to \kappa , $$
$$
k_i:\ z_i\to q^{-1/2}z_i,\ \ \
{\hat z}_i\to q^{1/2}{\hat z}_i,\ \ \
z_j\to z_j,\ \  {\hat z}_j\to {\hat z}_j, \ \ j\ne i-1,i.
$$
This action turns $M_q^{2n-1}$ into $U_q({\rm u}_n(\epsilon ))$-module
algebra. Harmonic analysis on the hyperboloid $M_q^{2n-1}$ demands to have
irreducible $*$-representations of $U_q({\rm u}_n(\epsilon ))$.
The harmonic analysis on the quantum hyperboloid $M^3_q$, related to the
quantum algebra $U_q({\rm su}_{1,1})$, see in [12].

In order to describe representations of the algebra $U_q({\rm u}_{n,1})$
we need finite dimensional representations of the subalgebra
$U_q({\rm u}_{n})$.
\medskip

\noindent {\bf 3. Irreducible representations of
$U_q({\rm u}_{n})$}
\medskip

The irreducible finite dimensional representations of
the algebra $U_q({\rm u}_{n})$
are given by $n$ integers ${\bf m}=(m_1,m_2,\cdots ,m_n)$ such that
$$
m_1\ge m_2\ge \cdots \ge m_n.
$$
These representations will be denoted by $T_{\bf m}$.
The set of numbers ${\bf m}=(m_1,m_2$, $\cdots$, $m_n)$ is called
the highest weight of the representation $T_{\bf m}$.
The representations $T_{\bf m}$ and $T_{{\bf m}'}$ are not equivalent if
${\bf m}\ne {\bf m}'$.

The Gel'fand--Tsetlin bases of carrier spaces of
irreducible representations $T_{\bf m}$
are formed by successive
restrictions of the representations to the subalgebras
$U_q({\rm u}_{n-1})$, $U_q({\rm u}_{n-2})$,
$\cdots ,U_q({\rm u}_1)\equiv U({\rm u}_1)$.
The decomposition of the
representation $T_{\bf m}$ of $U_q({\rm u}_n)$
into irreducible representations of $U_q({\rm u}_{n-1})$ is the same as
for the corresponding
representation $T_{\bf m}$ of ${\rm gl}(n,{\bf C})$.
Hence the restriction of $T_{\bf m}$,
${\bf m}=(m_1,\cdots ,m_n)$, to $U_q({\rm u}_{n-1})$ decomposes
into the irreducible representations $T_{{\bf m}_{n-1}}$,
${\bf m}_{n-1}=(m_{1,n-1},\cdots$, $m_{n-1,n-1})$,
such that
$$
m_1\ge m_{1,n-1}\ge m_2\ge m_{2,n-1}\ge \cdots \ge m_{n-1,n-1}
\ge m_n \eqno (9)
$$
and each of these representations enters into the decomposition
exactly once. Since the irreducible representations of $U({\rm u}_1)$ are
one-dimensional, we obtain a basis of the carrier
space $V_{\bf m}$ of the representation $T_{\bf m}$ of
$U_q({\rm u}_n)$ labeled by the Gel'fand--Tsetlin tableaux
$$
M=\left( \matrix{
      m_{1,n}\! \! & & m_{2,n} & & \cdots &  \cdots   & &\! \! m_{n,n}\cr
       & m_{1,n-1} & & m_{2,n-1} &  \cdots & & m_{n-1,n-1}\cr
         & & \cdots & \cdots  &\cdots  & \cdots & &\cr
         & & & & m_{11} &  & & &} \right) ,
\eqno (10)
$$
where $m_{i,n}\equiv m_i$.
The entries in (10) are integers satisfying the betweenness conditions
$$
m_{i,j+1}\ge m_{ij}\ge m_{i+1,j+1} ,\qquad i=1,2,\cdots ,j,\quad
j=1,2,\cdots ,n-1. \eqno (11)
$$
The set of all tableaux (10), satisfying these conditions,
labels the basis elements of the carrier space of
$T_{\bf m}$. If $M$ is a tableau (10), then the corresponding basis element
will be denoted by $| M\rangle $.

It was stated by Jimbo [10] that
the generators of $U_q({\rm u}_n)$ act on
the Gel'fand--Tsetlin basis of the representation by the formulas
$$
T_{\bf m}(k_r)| M \rangle =q^{a_r} | M\rangle ,\quad
a_r=\sum _{i=1}^r m_{ir}-\sum _{i=1}^{r-1} m_{i,r-1},\quad
1\le r\le n,\eqno (12) $$
$$
T_{\bf m}(e_r) | M\rangle =\sum _{j=1}^r A^j_r(M)| M^{+j}_r\rangle ,\ \
T_{\bf m}(f_r) | M\rangle =\sum _{j=1}^r A^j_r(M_r^{-j})
| M^{-j}_r\rangle ,
\eqno (13)$$
$$
\qquad\qquad\qquad\qquad\qquad\qquad\qquad  1\le r\le n-1.
$$
Here $M^{\pm j}_r$ is the Gel'fand--Tsetlin tableau obtained from the
tableau (10) if $m_{jr}$ is replaced by $m_{jr}\pm 1$, and
$A^j_r(M)$ is the expression
$$
A^j_r(M) =\left( - \frac {\prod _{i=1}^{r+1} [l_{i,r+1}-l_{jr}]
\prod _{i=1}^{r-1} [l_{i,r-1}-l_{jr}-1]}
{\prod _{i\ne j} [l_{ir}-l_{jr}][l_{ir}-l_{jr}-1]}
\right) ^{1/2} , \eqno (14)
$$
where $l_{is}=m_{is}-i$,
the positive value of the square root is taken and a number in square
brackets denotes a $q$-number defined by
$$
[m]=\frac {q^m-q^{-m}}{q-q^{-1}}.
$$

\noindent {\bf 4. Definition of representations of
$U_q({\rm u}_{n,1})$}

\medskip
In this paper
we are interested in representations of the quantum algebra
$U_q({\rm u}_{n,1})$. We destinguish for this algebra the notions of
representations and $*$-representations. Roughly speaking, the
notion of a representation does not take into account the $*$-structure of
$U_q({\rm u}_{n,1})$. A $*$-representation is a representation conserving
the $*$-structure. The strict definitions are as follows.

    A representation $T$ of the algebra $U_q({\rm u}_{n,1})$
is an algebraic homomorphism from $U_q({\rm u}_{n,1})$ to an algebra of
linear (bounded or unbounded) operators on a Hilbert space ${\cal H}$
for which the following conditions are fulfilled:
\medskip

    (a) the restriction of $T$ onto the maximal compact
subalgebra $U_q({\rm u}_n+{\rm u}_1)$
decomposes into a direct sum of its finite dimensional irreducible
representations;
\smallskip

    (b) operators of a representation $T$ are defined on everywhere dense
subspace ${\cal V}$ of ${\cal H}$, containing all subspaces which are carrier spaces of
irreducible finite dimensional subrepresentations of
$U_q({\rm u}_n+{\rm u}_1)$
from the restriction of $T$ onto this subalgebra;
\smallskip

(c) the subspace ${\cal V}$ is invariant with respect to all operators of
the representation $T$.
\medskip

To determine a representation $T$ of $U_q({\rm u}_{n,1})$
it is sufficient to give the operators $T(k_i^{\pm 1})$, $T(e_j)$, $T(f_j)$,
$i=1,2,\cdots ,n+1$, $j=1,2,\cdots ,n$, satisfying the relations (1)--(6)
written down for the algebra $U_q({\rm gl}_{n+1})$.

A representation $T$ of $U_q({\rm u}_{n,1})$ is called
$U_q({\rm u}_n+{\rm u}_1)$-finite if each irreducible representation
of $U_q({\rm u}_n+{\rm u}_1)$ is contained in the restriction of $T$ to
$U_q({\rm u}_n+{\rm u}_1)$ with finite multiplicity.
In other words, $U_q({\rm u}_n+{\rm u}_1)$-finite representations
of $U_q({\rm u}_{n,1})$ are
Harish-Chandra modules of $U_q({\rm u}_{n,1})$ with respect to
$U_q({\rm u}_n+{\rm u}_1)$. A linear span of irreducible subspaces of
${\cal H}$ with respect to the subalgebra $U_q({\rm u}_n+{\rm u}_1)$
will be denoted by ${\cal D}$. We shall see that ${\cal D}$
is invariant with respect to
$U_q({\rm u}_n+{\rm u}_1)$-finite representations of
$U_q({\rm u}_{n,1})$.

Below we consider only
$U_q({\rm u}_n+{\rm u}_1)$-finite representation of
$U_q({\rm u}_{n,1})$ and shall not emphasize this every time. So, everywhere
below under a representation of $U_q({\rm u}_{n,1})$ we understand its
$U_q({\rm u}_n+{\rm u}_1)$-finite representation.

A representation $T$ of $U_q({\rm u}_{n,1})$ on ${\cal H}$ is called
irreducible if the subspace ${\cal D}$ has no non-trivial invariant
subspaces (that is, ${\cal H}$ has no non-trivial invariant subspaces
with closure coinciding with ${\cal H}$).

If operators of a representation $T$ satisfy relations (7) and (8) on
the common domain of definition ${\cal V}$, then $T$ is called a
$*$-representation.

Two representations $T$ and $T'$ of $U_q({\rm u}_{n,1})$ on Hilbert spaces
${\cal H}$ and ${\cal H}'$ are called (algebraically) equivalent if there
exists a one-to-one operator $A: {\cal D}\to {\cal D}'$ such that
$T'(a)Av=AT(a)v$ for all $a\in U_q({\rm u}_{n,1})$ and $v\in {\cal D}$.
\medskip

\noindent
{\bf Lemma 1.} {\it Let $T$ be a representation of $U_q({\rm u}_{n,1})$
on a Hilbert space ${\cal H}$ such that the restriction of $T$ to the
subalgebra
$U_q({\rm u}_{n})$ decomposes into a direct sum of its irreducible
finite dimensional representations $T_{{\bf m}_n}$
and each of them is contained in the
decomposition with the unit multiplicity. Let
$$
{\cal H}=\oplus _{{\bf m}_n} {\cal V}_{{\bf m}_n} \eqno (15)
$$
be the corresponding decomposition of ${\cal H}$
into irreducible $U_q({\rm u}_{n})$-invariant subspaces. If
${\cal H}'$ is a $U_q({\rm u}_{n,1})$-invariant subspace of ${\cal H}$,
then ${\cal H}'$ is a direct sum of some subspaces ${\cal V}_{{\bf m}_n}$
of the decomposition (15).}
\medskip

\noindent
{\sl Proof.} Let ${\cal H}'$ be a $U_q({\rm u}_{n,1})$-invariant
subspace of ${\cal H}$ and let $v=v_{{\bf m}_n}+v_{{\bf m}'_n}\in {\cal H}'$,
where $v_{{\bf m}_n}\in {\cal V}_{{\bf m}_n}$ and
$v_{{\bf m}'_n}\in {\cal V}_{{\bf m}'_n}$ (for simplicity, we consider the
case of two summands; the case of more summands is considered similarly).
By the results of [13] (see also Proposition 7.9 in [5]), there exists an
element $Z$ of the center of $U_q({\rm u}_{n})$ such that
$\lambda _{{\bf m}_n}\equiv T_{{\bf m}_n}(Z)\ne
\lambda _{{\bf m}'_n}\equiv T_{{\bf m}'_n}(Z)$. Since
$$
T(Z)v=T_{{\bf m}_n}(Z)v_{{\bf m}_n}+
T_{{\bf m}'_n}(Z)v_{{\bf m}'_n}=
\lambda _{{\bf m}_n}v_{{\bf m}_n}+
\lambda _{{\bf m}'_n}v_{{\bf m}'_n}\in {\cal H}',
$$
then $T(Z)v-\lambda _{{\bf m}_n}v=
(\lambda _{{\bf m}'_n}-\lambda _{{\bf m}_n})v_{{\bf m}'_n}\in {\cal H}'$.
Thus, $v_{{\bf m}_n}\in {\cal H}'$ and $v_{{\bf m}'_n}\in {\cal H}'$.
This means that $T(U_q({\rm u}_{n})) v_{{\bf m}_n}={\cal V}_{{\bf m}_n}
\in {\cal H}'$ and
$T(U_q({\rm u}_{n})) v_{{\bf m}'_n}={\cal V}_{{\bf m}'_n}
\in {\cal H}'$. Lemma is proved.
\medskip

{\bf 5. The principal series representations of $U_q({\rm u}_{n,1})$}
\medskip

In this section we construct a series of representations of
$U_q({\rm u}_{n,1})$ which is a $q$-analogue of the principal nonunitary
series of the real Lie algebra ${\rm u}_{n,1}$. We first describe these
representations and then prove that they indeed are representations of
$U_q({\rm u}_{n,1})$.

Let
$c_1$ and $c_2$ be complex numbers such that $c_1+c_2=m_0$ is an integer,
and let ${\bf m}=(m_1,m_2,\cdots ,m_{n-1})$ be a set if integers such that
$m_1\ge m_2\ge \cdots \ge m_{n-1}$. The number $m_0$ determines a
one-dimensional representation of the subalgebra $U_q({\rm u}_{1})$ and
${\bf m}$ determines an irreducible finite dimensional representation
of the subalgebra $U_q({\rm u}_{n-1})$.
The numbers ${\bf m}$, $c_1$, $c_2$ determine the representation
$T_{{\bf m},c_1,c_2}$ of
the quantum algebra $U_q({\rm u}_{n,1})$ which is defined as follows.
The restriction of $T_{{\bf m},c_1,c_2}$ onto $U_q({\rm u}_n)$
decomposes into the direct sum of all irreducible representations
$T_{{\bf m}_n}$ of this subalgebra with highest weights
${\bf m}_n=(m_{1n},\cdots ,m_{nn})$ such that
$$m_{1n}\ge m_1\ge m_{2n}\ge m_2\ge \cdots \ge m_{n-1}\ge m_{nn}.\eqno (16)$$
Every of these representations of $U_q({\rm u}_n)$ is contained in
$T_{{\bf m},c_1,c_2}$ exactly once. This restriction determines the
Hilbert space ${\cal H}_{\bf m}$ (independent of $c_1$ and $c_2$)
on which $T_{{\bf m},c_1,c_2}$ acts. This space is described by the
decomposition ${\cal H}_{\bf m} =\oplus _{{\bf m}_n} {\cal V}_{{\bf m}_n}$,
where ${\cal V}_{{\bf m}_n}$ is the subspace on which the irreducible
representation of $U_q({\rm u}_n)$ with highest weight ${\bf m}_n$
is realized and the sum is over all highest weights satisfying (16).
We choose in ${\cal H}_{\bf m}$ the orthonormal basis consisting of the
Gel'fand-Tsetlin bases of the subspaces ${\cal V}_{{\bf m}_n}$.
The basis elements are denoted by $\vert M\rangle \equiv \vert {\bf m}_n$,
$\alpha \rangle$, where $M$ is a Gel'fand--Tsetlin tableau (10) and
$\alpha$ is the tableau $M$ without the first row ${\bf m}_n$.
The operators
$T_{{\bf m},c_1,c_2}(e_n)$ and $T_{{\bf m},c_1,c_2}(f_n)$ of the
representation $T_{{\bf m},c_1,c_2}$ act upon
$\vert {\bf m}_n,\alpha \rangle$ by the formulas
$$T_{{\bf m},c_1,c_2}(e_n)\vert {\bf m}_n,\alpha \rangle
=\sum ^n_{s=1} [l_{sn}-c_1]
\omega _s({\bf m},{\bf m}_n,\alpha )\vert {\bf m}_n^{+s},\alpha
\rangle , \eqno (17)$$
$$T_{{\bf m},c_1,c_2}(f_n)\vert {\bf m}_n,\alpha \rangle
=\sum ^n_{s=1} [-l_{sn}+c_2+1]
\omega _s({\bf m},{\bf m}_n^{-s},\alpha )\vert {\bf m}_n^{-s},
\alpha \rangle , \eqno (18)$$
where
$$\omega _s({\bf m},{\bf m}_n,\alpha )=\Biggl( {\prod _{j=1}^{n-1}
[l_{j,n-1}-l_{sn}-1][l_j-l_{sn}]\over \prod _{r\ne s} [l_{sn}-l_{rn}+1]
[l_{sn}-l_{rn}]}\Biggr) ^{1/2}, \eqno (19)$$
$$l_j=m_j-j-1,\ \ j=1,2,\cdots ,n-1;\ \ \ l_{sk}=m_{sk}-s,\ \ s=1,2,
\cdots ,k,$$
and ${\bf m}_n^{\pm s}=(m_{1n},\cdots m_{s-1,n},m_{sn}\pm 1, m_{s+1,n},
\cdots ,m_{nn})$ if ${\bf m}_n=(m_{1n},\cdots ,m_{nn})$.
The operator $T_{{\bf m},c_1,c_2}(k_{n+1})$ is given by the formula
$$
T_{{\bf m},c_1,c_2}(k_{n+1})\vert {\bf m}_n,\alpha \rangle =q^a
\vert {\bf m}_n,\alpha \rangle ,\ \  \
a =\Bigl( c_1+c_2+n+2
+\sum _{j=1}^{n-1}m_j-\sum _{j=1}^n m_{jn} \Bigr) . \eqno (20)
$$
The other generators $e_i,\ f_j,\ h_k$ belong to the subalgebra
$U_q({\rm u}_n)$ and the operators $T_{{\bf m},c_1,c_2}(e_i)$,
$T_{{\bf m},c_1,c_2}(f_j)$, $T_{{\bf m},c_1,c_2}(h_k)$ act upon
the basis elements $\vert {\bf m}_n,\alpha \rangle$ by the corresponding
formulas (12)--(14).
\medskip

\noindent
{\bf Theorem 1.} {\it The mapping $a\to T_{{\bf m},c_1,c_2} (a)$,
$a\in U_q({\rm u}_{n,1})$, described above, indeed determines a
representation of the algebra $U_q({\rm u}_{n,1})$.}
\medskip

\noindent
{\sl Proof.} For simplicity, in this proof we denote the operators
$T_{{\bf m},c_1,c_2}(e_i)$, $T_{{\bf m},c_1,c_2}(f_j)$ and
$T_{{\bf m},c_1,c_2}(k_r)$ by $e_i$, $f_j$ and $k_r$, respectively.
We have to prove the following relations ($1 \leq i \leq n,\ $
$1 \leq j \leq n-2$) :
$$
k_{n+1} k_{n+1}^{-1}= k_{n+1}^{-1}k_{n+1}=1,\ \ \
k_{n+1} k_i = k_i k_{n+1},      \eqno (21)
$$
$$
k_{n+1} e_i  k_{n+1}^{-1} = q^{-\delta_{i,n}} e_i,\ \
k_{n+1} f_i  k_{n+1}^{-1} = q^{\delta_{i,n}} f_i,\ \
k_{i} e_n  k_{i}^{-1} = q^{\delta_{i,n}} e_n,\ \
k_{i} f_n  k_{i}^{-1} = q^{-\delta_{i,n}} f_n, \eqno (22)
$$
$$
[e_n,f_{n-1}]=[e_{n-1},f_n]=
[e_n,e_j]=[f_n,f_j]=[e_n,f_j]=[e_j,f_n]=0, \eqno (23)
$$
$$
e_n^2e_{n-1}+e_{n-1}e_n^2-[2] e_n e_{n-1}e_n =0,\ \ \
e_{n-1}^2 e_n+e_n e_{n-1}^2-[2] e_{n-1} e_{n} e_{n-1} =0, \eqno (24)
$$
$$
f_n^2f_{n-1}+f_{n-1}f_n^2-[2] f_n f_{n-1}f_n =0, \ \ \
f_{n-1}^2 f_n+f_n f_{n-1}^2-[2] f_{n-1} f_{n} f_{n-1} =0, \eqno (25)
$$
$$
[e_n,f_n]=\frac{k_n k_{n+1}^{-1} - k_{n+1}k_n^{-1}}
{q-q^{-1}}. \eqno (26)
$$
The relations (21) are trivial. The relations (22) are easily verified
by means of formulas (12)--(14) and
(17)--(20). The relations (23) are also easily verified.
As a sample, we check the relation $[e_n,f_{n-1}]=0$.
Acting by the left hand side upon the basis vector $\vert M\rangle :=
\vert {\bf m}_n,\alpha\rangle$ and then collecting all coefficients
at the basis vector $|(M_n^{+j'})_{n-1}^{-j}\rangle$ we obtain the
expression
$$
[l_{j',n}-c_1]\omega_{j'}(M_{n-1}^{-j})
A_{n-1}^j(M_{n-1}^{-j})-
A_{n-1}^j((M_n^{+j'})_{n-1}^{-j})[l_{j',n}-c_1]
\omega_{j'}(M) .
$$
Substituting here the expressions for $\omega _{j'}$ and $A^j_{n-1}$
we reduce it to the expression
$$
D([l_{j,n-1}-l_{j',n}-2][l_{j',n}-l_{j,n-1}+1]-
[l_{j,n-1}-l_{j',n}-1][l_{j',n}-l_{j,n-1}+2]),
$$
where $D$ stands for a certain expression, explicit form of which is not
important for us. The expression in the brackets trivially vanishes.

All relations in (24) and (25) are verified in the same manner. As a sample,
we check the relation
$$
e_n^2e_{n-1}+e_{n-1}e_n^2-[2] e_n e_{n-1}e_n =0. \eqno (27)
$$
We act by the left hand side upon the basis vector $\vert M\rangle$.
As a result we obtain a linear combinations of basis vectors of
the types $|((M_n^{+j'})_n^{+j'})_{n-1}^{+j}\rangle$ and
$|((M_n^{+j'})_n^{+j''})_{n-1}^{+j}\rangle,\ j' \neq j''$.
Collecting the coefficients at the vector
$|((M_n^{+j'})_n^{+j'})_{n-1}^{+j}\rangle$ we obtain the expression
$$
[(l_{j',n}+1)-c_1]\omega_{j'}((M_n^{+j'})_{n-1}^{+j})
[l_{j',n}-c_1]\omega_{j'}(M_{n-1}^{+j})A_{n-1}^j(M)+
$$
$$
+A_{n-1}^j((M_n^{+j'})_{n}^{+j'})[(l_{j',n}+1)-c_1]
\omega_{j'}(M_n^{+j'})
[l_{j',n}-c_1]\omega_{j'}(M)-
$$
$$
-[2][(l_{j',n}+1)-c_1]\omega_{j'}((M_n^{+j'})_{n-1}^{+j})
A_{n-1}^j(M_{n}^{+j'})[l_{j',n}-c_1]\omega_{j'}(M). \eqno (28)
$$
Substituting the expressions for $\omega _{j'}$ and $A^j_{n-1}$
we reduce it (up to some common coefficient) to the expression
$$
(-[l_{j,n-1}-l_{j',n}-1][l_{j,n-1}-l_{j',n}]
[l_{j',n}-l_{j,n-1}])^{1/2}+
$$
$$
+(-[l_{j,n-1}-l_{j',n}-2][l_{j,n-1}-l_{j',n}-1]
[l_{j',n}-l_{j,n-1}+2])^{1/2}-
$$
$$
-[2](-[l_{j,n-1}-l_{j',n}-1][l_{j,n-1}-l_{j',n}-1]
[l_{j',n}-l_{j,n-1}+1])^{1/2} .
$$
It reduces to
$$
([l_{j,n-1}-l_{j'n}-1])^{1/2}([x+2]+[x]-[2][x+1])\ \ \ {\rm where}\ \ \
x=l_{j,n-1}-l_{j'n}-2.
$$
It is easy to check that
$$
[x+2]+[x]-[2][x+1]=0.  \eqno (29)
$$
That is, the expression (28) vanishes.

Collecting the coefficients at the basis vector
$|((M_n^{+j'})_n^{+j''})_{n-1}^{+j}\rangle,\ j' \neq j'',$
and (for a convenience) dividing
them by $A_{n-1}^j$ we obtain the expression
$$
\omega_{j'}(M_n^{+j''})\omega_{j''}(M)
\left(
\frac{\omega_{j'}((M_n^{+j''})_{n-1}^{+j})}
{\omega_{j'}(M_n^{+j''})}
\frac{\omega_{j''}(M_{n-1}^{+j})}
{\omega_{j''}(M)}+
\frac{A_{n-1}^j((M_n^{+j'})_{n}^{+j''})}
{A_{n-1}^j(M)}\right.
$$
$$
\left.
-[2]\frac{\omega_{j'}((M_n^{+j''})_{n-1}^{+j})}
{\omega_{j'}(M_n^{+j''})}
\frac{A_{n-1}^j(M_n^{+j''})}
{A_{n-1}^j(M)}\right)+
\omega_{j''}(M_n^{+j'})\omega_{j'}(M)
\Biggl(\ \  j'\leftrightarrow  j''\ \  \Biggr) ,
                                                  \eqno(30)
$$
where $(\ \  j'\leftrightarrow  j''\ \  )$ means the expression in
the first big round brackets with $j'$ replaced by $j''$ and $j''$ by
$j'$. Taking the explicit forms for
$\omega_s$ and $A_{n-1}^j$, and then using the identity (29) and
the equalities
$$
[x][y]=[(x+y)/2]^2-[(x-y)/2]^2,\ \ \
[x]^2-[y]^2=[x+y][x-y]
$$
we reduce the expression in
the first big round brackets to
$$
-[l_{j',n}-l_{j'',n}-1]\bigl([l_{j,n-1}-l_{j',n}]
[l_{j,n-1}-l_{j'',n}][l_{j,n-1}-l_{j',n}-1]
[l_{j,n-1}-l_{j'',n}-1]\bigr)^{-1/2}.    \eqno(31)
$$
The expression in the second big round brackets can
be obtained from (31) by replacement $j'\leftrightarrow  j''$.
Substituting these expressions and the explicit forms for
$\omega_{j'}$ and $\omega_{j''}$ into (30) we represent (30) in
the form
$$
N\left( \frac{[l_{j',n}-l_{j'',n}-1]}
{([l_{j',n}-l_{j'',n}][l_{j',n}-l_{j'',n}-1]
[l_{j'',n}-l_{j',n}+1][l_{j'',n}-l_{j',n}])^{1/2}}+ \right.
$$
$$\left.
+\frac{[l_{j'',n}-l_{j',n}-1]}
{([l_{j'',n}-l_{j',n}][l_{j'',n}-l_{j',n}-1]
[l_{j',n}-l_{j'',n}+1][l_{j',n}-l_{j'',n}])^{1/2}} \right) , \eqno (32)
$$
where $N$ stands for a certain expression, explicit form of which is not
important for us.
It is easily checked that the expression in the brackets (at the coefficient
$N$) vanishes. Therefore the expression (30) vanishes and therefore
the relation (27) is proved.

Now let us prove the relation (26). Acting by the left hand side upon
the basis vector $\vert M\rangle$ we obtain a linear combination of the
basis vectors $|(M_n^{+j'})_n^{-j''}\rangle$, $j'\ne j''$, and
$\vert M\rangle$.
Acting by the right hand side of (26) upon
the basis vector $\vert M\rangle$ we obtain the same basis vector with
a coefficient.

Acting by both sides of (26) upon $\vert M\rangle$ and collecting
coefficients at the basis vector
$|(M_n^{+j'})_n^{-j''}\rangle$, $j'\ne j''$ we obtain the
relation
$$
\omega_{j'}(M_n^{-j''})\omega_{j''}(M_n^{-j''})
-\omega_{j''}((M_n^{+j'})_n^{-j''})\omega_{j'}(M)=0.
$$
Using the expicit form of $\omega_{j'}$ and $\omega_{j''}$
it reduces to the relation
$$
\bigl([(l_{j'',n}-1)-l_{j',n}][(l_{j'',n}-1)-l_{j',n}-1]
[l_{j',n}-(l_{j'',n}-1)][l_{j',n}-(l_{j'',n}-1)-1]\bigr)^{-1/2}
$$
$$
-\bigl([(l_{j',n}+1)-(l_{j'',n}-1)][(l_{j',n}+1)-(l_{j'',n}-1)-1]
[l_{j'',n}-l_{j',n}][l_{j'',n}-l_{j',n}-1]\bigr)^{-1/2}=0
$$
which is obviously fulfilled.

Equating diagonal matrix elements of both sides of (26)
we obtain the relation
$$
\Phi(\{l_{1,n+1},l_{1,n},l_{1,n-1}\};
\{l_{2,n+1},l_{2,n},l_{2,n-1}\};\cdots ;
\{l_{n,n+1},l_{n,n}\}; \{l_{n+1,n+1}\}):=
$$
$$
:= \sum_{j=1}^n\Bigl((A_n^j(M_n^{-j}))^2-
(A_n^j(M))^2\Bigr)
=\left[ 2\sum_{j=1}^n l_{j,n}-\sum_{j=1}^{n+1} l_{j,n+1}-
\sum_{j=1}^{n-1} l_{j,n-1}-1 \right] ,
                                                      \eqno(33)
$$
where the notations $l_{1,n+1}:=c_1$, $l_{n+1,n+1}:=c_2$ and
$l_{j+1,n+1}:=l_j\ (1\leq j \leq n-1)$ were used and
$(A_n^j(M))^2$ is the square of the expression (14) taken for $r=n$.
We consider the relation (33) for arbitrary complex
$l_{1,n+1}$ and $l_{n+1,n+1},$ with their sum not necessary
to be integer.

We first suppose that the relation (33) is valid for
$l_{1,n+1}=l_{1,n}$ and show that then (33) is true for any value of
$l_{1,n+1}$. For this, we derive some properties of the function $\Phi$.

Using the identity $[a][b]-[a+x][b-x]=[x][a-b+x],$
which is valid for arbitrary complex $a,b,x,$ it
is directly verified the following identity for
function $\Phi$:
$$
\Phi(\{l_{1,n+1},l_{1,n},l_{1,n-1}\};
\cdots ;\{l_{n,n+1},l_{n,n}\}; \{l_{n+1,n+1}\})-
$$
$$
-\Phi(\{l_{1,n+1}+x,l_{1,n},l_{1,n-1}\};
\cdots ;\{l_{n,n+1},l_{n,n}\}; \{l_{n+1,n+1}-x\})=
$$
$$
=[x][l_{1,n+1}-l_{n+1,n+1}+x]
\Phi(\{ \cdot ,l_{1,n},l_{1,n-1}\} ;
\cdots ;\{l_{n,n+1},l_{n,n}\}; \{ \cdot \}) ,  \eqno(34)
$$
where the dot "$\cdot $" on the places of $l_{1,n+1}$ and $l_{n+1,n+1}$
in the last $\Phi$ means that all the multipliers depending on
$l_{1,n+1}$ and $l_{n+1,n+1}$ in the expression for $\Phi$ must be omitted.

Now let us prove the identity
$$
\Phi(\{l_{1,n}-1,l_{1,n},l_{1,n-1}\};\cdots )=
\Phi(\{l_{1,n}+1,l_{1,n}+1,l_{1,n-1}\}; \cdots )
                                            \eqno(35)
$$
For summands in $\Phi(\{l_{1,n}-1,l_{1,n},l_{1,n-1}\} ;\cdots )$ we have
$$
(A_n^1(\{l_{1,n}-1,l_{1,n}-1,l_{1,n-1}\}; \cdots ))^2=0 ,
$$
$$
(A_n^1(\{l_{1,n}-1,l_{1,n},l_{1,n-1}\}; \cdots ))^2=
$$
$$
=-\frac{[(l_{1,n}-1)-l_{1,n}]\prod_{s=2}^{n+1}[l_{s,n+1}-l_{1,n}]
\prod_{s=1}^{n-1}[l_{s,n-1}-l_{1,n}-1]}
{\prod_{s=2}^{n}[l_{s,n}-l_{1,n}][l_{s,n}-l_{1,n}-1]} ,
$$
$$
(A_n^j(\{l_{1,n}-1,l_{1,n},l_{1,n-1}\}; \cdots ))^2=
$$
$$
=-\frac{[(l_{1,n}-1)-l_{j,n}]\prod_{s=2}^{n+1}[l_{s,n+1}-l_{j,n}]
\prod_{s=1}^{n-1}[l_{s,n-1}-l_{j,n}-1]}
{[l_{1,n}-l_{j,n}][l_{1,n}-l_{j,n}-1]
\prod_{s\geq 2,s\neq j}[l_{s,n}-l_{j,n}][l_{s,n}-l_{j,n}-1]},
\ \ \ 2\ge j \ge n ,
$$
and for summands in $\Phi(\{l_{1,n}+1,l_{1,n}+1,l_{1,n-1}\}; \cdots )$
we obtain
$$
(A_n^1(\{l_{1,n}+1,l_{1,n},l_{1,n-1}\}; \cdots ))^2=
$$
$$
=-\frac{\prod_{s=2}^{n+1}[l_{s,n+1}-l_{1,n}]
\prod_{s=1}^{n-1}[l_{s,n-1}-l_{1,n}-1]}
{\prod_{s=2}^{n}[l_{s,n}-l_{1,n}][l_{s,n}-l_{1,n}-1]} ,
$$
$$
(A_n^1(\{l_{1,n}+1,l_{1,n}+1,l_{1,n-1}\}; \cdots ))^2=0,
$$
$$
(A_n^j(\{l_{1,n}+1,l_{1,n}+1,l_{1,n-1}\}; \cdots ))^2=
$$
$$
=-\frac{\prod_{s=2}^{n+1}[l_{s,n+1}-l_{j,n}]
\prod_{s=1}^{n-1}[l_{s,n-1}-l_{j,n}-1]}
{[l_{1,n}-l_{j,n}]
\prod_{s\ge 2,s\ne j}[l_{s,n}-l_{j,n}][l_{s,n}-l_{j,n}-1]}
,\ \ \ 2\ge j \ge n .
$$
We easily see that the above expressions for
the summands in $\Phi(\{l_{1,n}-1,l_{1,n},l_{1,n-1}\}; \cdots )$
coincide with those for
the summands in $\Phi(\{l_{1,n}+1,l_{1,n}+1,l_{1,n-1}\}; \cdots )$.
Thus, the relation (35) is proved.

From the identity (35) with
$l_{n+1,n+1}$ replaced by $l_{n+1,n+1}+1$ and
from the assumption that (33) is true for $l_{1,n+1}=l_{1,n}$
we derive that
$$
\Phi(\{l_{1,n}-1,l_{1,n},l_{1,n-1}\} ;
\cdots ;\{l_{n,n+1},l_{n,n}\}; \{l_{n+1,n+1}+1\}) =
$$
$$
=\Phi(\{l_{1,n}+1,l_{1,n}+1,l_{1,n-1}\} ;
\cdots ;\{l_{n,n+1},l_{n,n}\}; \{l_{n+1,n+1}+1\}) =
$$
$$
=\left[ l_{1,n}+2\sum_{j=2}^n l_{j,n}
-\sum_{j=2}^{n+1} l_{j,n+1}
-\sum_{j=1}^{n-1} l_{j,n-1}-1 \right] =
$$
$$
=\Phi(\{l_{1,n},l_{1,n},l_{1,n-1}\} ;
\cdots ;\{l_{n,n+1},l_{n,n}\}; \{l_{n+1,n+1}\})  .
$$
It follows from here that at
$l_{1,n+1}=l_{1,n}$ and $x=-1$ ($l_{n+1,n+1}$ is an
arbitrary complex number) the left hand side of (34) identically vanishes.
Hence the right hand side vanishes. Since
$l_{n+1,n+1}$ is an arbitrary complex number, we obtain that
$\Phi$ from the right hand side vanishes:
$$
\Phi(\{ \cdot ,l_{1,n},l_{1,n-1}\} ;
\cdots ;\{l_{n,n+1},l_{n,n}\}; \{ \cdot \}) = 0.
$$
Using again the identity (34)
with arbitrary $l_{1,n+1}$ and with $x=l_{1,n}-l_{1,n+1}$ and then
the relation (33) at $l_{1,n+1}=l_{1n}$ we have
$$
\Phi(\{l_{1,n+1},l_{1,n},l_{1,n-1}\} ;
\cdots ;\{l_{n,n+1},l_{n,n}\}; \{l_{n+1,n+1}\})
$$
$$
=\Phi(\{l_{1,n},l_{1,n},l_{1,n-1}\} ;
\cdots ;\{l_{n,n+1},l_{n,n}\}; \{l_{n+1,n+1}+(l_{1,n+1}-l_{1,n})\})
$$
$$
=\left[ 2\sum_{j=1}^n l_{j,n}
-\sum_{j=1}^{n+1} l_{j,n+1}
-\sum_{j=1}^{n-1} l_{j,n-1}-1 \right] .
$$
This gives the relation (33) for arbitrary $l_{1,n+1}$.

It is shown similarly that the relation (33)
with the admissable values of $l_{1,n-1}$ follows from (33) with
$l_{1,n-1}=l_{1,n}$.

Thus, in order to prove the relation (33) with arbitrary $l_{1,n+1}$
$l_{1,n},$ $l_{1,n-1}$ it is enough to prove the special
case of (33) when $l_{1,n+1}=l_{1,n}=l_{1,n-1}$. In order to
prove this special case we note that
$$
\Phi(\{l_{1,n},l_{1,n},l_{1,n}\} ;
\{l_{2,n+1},l_{2,n},l_{2,n-1}\}; \cdots )=
\Phi(\{l_{2,n+1},l_{2,n},l_{2,n-1}\}; \cdots ),         \eqno(36)
$$
where on the right hand side
all the multipliers depending on $l_{1,n+1},$
$l_{1,n},$ $l_{1,n-1}$ are omitted.
Indeed, we have the equalities
$$
(A_n^1(\{l_{1,n},l_{1,n}-1,l_{1,n}\};\cdots ))^2=0 ,\
\ \ (A_n^1(\{l_{1,n},l_{1,n},l_{1,n}\};\cdots ))^2=0 ,
$$
$$
(A_n^j(\{l_{1,n},l_{1,n},l_{1,n}\} ;
\{l_{2,n+1},l_{2,n},l_{2,n-1}\};\cdots ))^2=
(A_n^j(\{l_{2,n+1},l_{2,n},l_{2,n-1}\};\cdots ))^2 ,
$$
with $2 \le j \le n$, which prove (36).

The right hand side of (33) at $l_{1,n+1}=l_{1,n}=l_{1,n-1}$
coincides with $[2\sum_{j=2}^n l_{j,n}-\sum_{j=2}^{n+1} l_{j,n+1}-
\sum_{j=2}^{n-1} l_{j,n-1}-1]$. This means that
the relation (33) at $l_{1,n+1}=l_{1,n}=l_{1,n-1}$
coincides with the same relation written down for the algebra
$U_q({\rm u}_{n-1,1})$.
Hence, the validity of the relation
(33), which correspond to $U_q({\rm u}_{n,1})$, follows from the validity
of (33) for $U_q({\rm u}_{n-1,1})$. If we continue the described procedure
we reduce the proof of (33) to the obvious
identity
$$
\Phi(\{l_{n,n},l_{n,n}\}; \{l_{n+1,n+1}\}) \equiv
- [l_{n,n}-l_{n,n}+1][l_{n+1,n+1}-l_{n,n}+1]=
[2l_{n,n}-l_{n,n}-l_{n+1,n+1}-1].
$$
Theorem is proved.
\medskip

There exist equivalence relations in the set of the representations
$T_{{\bf m},c_1,c_2}$. A part of these equivalences is related to
the periodicity of analytical function
$$
w(z)=[z]\equiv (q^z-q^{-z})/(q-q^{-1}).
$$
We set $q=\exp h$. Then $w(z)$ is a periodic function with period
$2\pi {\rm i}/h$. Besides, we have $w(z)=-w(z+\pi {\rm i}/h)$. This means
that the representations
$$
T_{{\bf m},c_1,c_2} \ \ \ {\rm and}\ \ \ T_{{\bf m},c_1+2\pi {\rm i}k/h,
c_2-2\pi {\rm i}k/h},\ \ \ k\in {\bf Z},
$$
coincide and the representations
$$
T_{{\bf m},c_1,c_2}\ \ \ {\rm and}\ \ \ T_{{\bf m},c_1+\pi {\rm i}k/h,
c_2-\pi {\rm i}k/h},\ \ \ k\in {\bf Z}, \eqno (37)
$$
are equivalent. Therefore, we may restrict ourselves by consideration of
the representations $T_{{\bf m},c_1,c_2}$ with
$-\pi /2h < {\rm Im}\, c_1 \le \pi /2h$. However, it will be often
convenient to assume that ${\rm Im}\, c_1$ is defined modulo
$\pi /h$ and $c_2=m_0-c_1$.

In order to study irreducibility of the representations $T_{{\bf m},c_1,c_2}$
we prove the following assertion:
\medskip

\noindent
{\bf Proposition 1.} {\it The representation $T_{{\bf m},c_1,c_2}$
with $-\pi /2h < {\rm Im}\, c_1 \le \pi /2h$ is
irreducible if for each irreducible representation $T_{{\bf m}_n}$ of
$U_q({\rm u}_n)$, contained in the restriction of $T_{{\bf m},c_1,c_2}$
to $U_q({\rm u}_n)$, no of the numbers $l_{sn}-c_1$ and
$l_{sn}-c_2-1$, $s=1,2,\cdots ,n$, vanishes. If $l_{sn}-c_1=0$ for
some representation $T_{{\bf m}_n}$ such that $T_{{\bf m}^{+s}_n}$
is not contained in the restriction $T_{{\bf m},c_1,c_2}{\downarrow}
U_q({\rm u}_n)$
or/and $l_{sn}-c_2-1=0$  for some representation $T_{{\bf m}_n}$ such that
$T_{{\bf m}^{-s}_n}$ is not contained in the restriction
$T_{{\bf m},c_1,c_2}{\downarrow}U_q({\rm u}_n)$, then $T_{{\bf m},c_1,c_2}$
is also irreducible.}
\medskip

\noindent
{\sl Proof.} Since $q$ is a positive number, the number $l_{sn}-c_1$ (resp.
the number $l_{sn}-c_2-1$) vanishes if and only if the $q$-number
$[l_{sn}-c_1]$ (resp. the $q$-number $[l_{sn}-c_2-1]$) vanishes.
These $q$-numbers are parts of the coefficients at the vectors on
the right hand sides of (17) and (18). Multipliers
$\omega _s({\bf m},{\bf m}_n,\alpha )$ and
$\omega _s({\bf m},{\bf m}^{-s}_n,\alpha )$ in these coefficients do
not vanish. Thus, under the conditions of our Proposition the
coefficients at $\vert {\bf m}^{+s}_n, \alpha \rangle$ and
$\vert {\bf m}^{-s}_n, \alpha \rangle$ in formulas (17) and (18)
do not vanish.

Let ${\cal H}'$ be a $U_q({\rm u}_{n,1})$-invariant subspace of the space
${\cal H}_{\bf m}$ of the representation $T_{{\bf m},c_1,c_2}$. By Lemma 1,
${\cal H}'=\oplus _{{\bf m}'_n} {\cal V}_{{\bf m}'_n}$. We take a fixed
vector $\vert {{\bf m}'_n}, \alpha\rangle$ from this subspace.
It follows from (17) and (18) that the operators $T_{{\bf m},c_1,c_2}(e_j)$,
acting upon the vectors $\vert {\bf m}_n,\alpha\rangle$,
increase the numbers $m_{ij}$ and
the operators $T_{{\bf m},c_1,c_2}(f_j)$ decrease these numbers. Since
the coefficients in (17) and (18) do not vanish, then acting
successively upon
$\vert {{\bf m}'_n}, \alpha\rangle$ by the operators
$T_{{\bf m},c_1,c_2}(e_i)$ and $T_{{\bf m},c_1,c_2}(f_i)$, $i=1,2,\cdots ,
n$, we can obtain a vector $v=\sum a_{{\bf m}_n,\alpha '}
\vert {{\bf m}_n}, \alpha '\rangle$ with non-vanishing coefficient at
$\vert {{\bf m}_n}, \alpha '\rangle$, where ${\bf m}_n$ is any fixed highest
weight satisfying the condition (16). Since ${\cal H}'$ is
$U_q({\rm u}_{n,1})$-invariant, we have $v\in {\cal H}'$ and therefore
$\vert {\bf m}_n,\alpha '\rangle \in {\cal H}'$ by Lemma 1. Hence
${\cal H}'={\cal H}$. Proposition is proved.
\medskip

\noindent
{\bf Corollary 1.} {\it
The representation $T_{{\bf m},c_1,c_2}$
with $-\pi /2h < {\rm Im}\, c_1 \le \pi /2h$ is
irreducible if $c_1$ and $c_2$ are not integers or if both $c_1$ and $c_2$
are integers coinciding
with some of the numbers $l_1,l_2,\cdots , l_{n-1}$, $l_i=m_i-i-1$.}
\medskip

\noindent
{\sl Proof.} For such values of $c_1$ and $c_2$ the numbers from
Proposition 1 do not vanish. This proves our Corollary.
\medskip

Below we show that the conditions of Corollary 1 are also sufficient
for irreducibility of $T_{{\bf m},c_1,c_2}$.
\medskip

Note that the algebra $U_q({\rm u}_{n,1})$ has one-dimensional
representations which are in fact representations of $U_q({\rm gl}_{n+1})$.
As in the case of finite dimensional representations, new infinite
dimensional representations of $U_q({\rm u}_{n,1})$ may be obtained by taking
tensor products of the representations $T_{{\bf m},c_1,c_2}$ with
one-dimensional ones. But these new representations are not essential
and we do not consider them below.

Below we shall need the following
\medskip

\noindent
{\bf Lemma 2.} {\it Let $T$ be a representation of $U_q({\rm u}_{n,1})$
such that $T{\downarrow} U_q({\rm u}_{n})$ contains irreducible
representations of $U_q({\rm u}_n)$ not more than once.  Let the operator
$T(k_{n+1})$ be diagonal in the Gel'fand--Tsetlin basis
$\{ \vert {\bf m}_n,\alpha\rangle \}$
of the representation space, and let $T(e_n)$ and $T(f_n)$ are
given by
$$
T(e_n) \vert {\bf m}_n,\alpha\rangle =\sum _{s=1}^n C_s(
{\bf m}_n,\alpha ) \vert {\bf m}^{+s}_n,\alpha\rangle ,\ \
T(f_n) \vert {\bf m}_n,\alpha\rangle =\sum _{s=1}^n B_s(
{\bf m}_n,\alpha ) \vert {\bf m}^{-s}_n,\alpha\rangle .
$$
If the coefficients $C_s({\bf m}_n,\alpha )$ and $B_s({\bf m}_n,\alpha )$
vanish only if ${\bf m}^{+s}_n$ and ${\bf m}^{-s}_n$ are not contained
in $T{\downarrow} U_q({\rm u}_{n})$, then the representation $T$ is
irreducible.}
\medskip

\noindent
{\sl Proof} is the same as that of Proposition 1 and we omit it.
\medskip

\noindent
{\bf 6. Intertwining operators}
\medskip

\noindent
In section 5 we found that the representations
$T_{{\bf m},c_1,c_2}$
and $T_{{\bf m},c_1+\pi {\rm i}k/h,c_2-\pi {\rm i}k/h}$, $k\in {\bf Z}$,
are equivalent. Let us find other equivalence relations in the set of
the irreducible representations $T_{{\bf m},c_1,c_2}$.

Assume that the representations $T_{{\bf m},c_1,c_2}$ and
$T_{{\bf m}',c'_1,c'_2}$ are equivalent. Then the restrictions
$T_{{\bf m},c_1,c_2}{\downarrow}U_q({\rm u}_n)$ and
$T_{{\bf m}',c'_1,c'_2}{\downarrow}U_q({\rm u}_n)$ are equivalent
representations of $U_q({\rm u}_n)$. This means that these restrictions
consist of the same irreducible representations of $U_q({\rm u}_n)$.
It follows from (16) that it is possible if and only if
${\bf m}={\bf m}'$.

Equivalence of the irreducible representations
$T_{{\bf m},c_1,c_2}$ and $T_{{\bf m},c'_1,c'_2}$ means that there exists an
invertible operator $A: {\cal D}\to {\cal D}$ such that
$$
AT_{{\bf m},c_1,c_2}(a)=T_{{\bf m},c'_1,c'_2}(a)A,\ \ \ \
a\in U_q({\rm u}_{n,1}), \eqno (38)
$$
on ${\cal D}$. Since
$T_{{\bf m},c_1,c_2}{\downarrow}U_q({\rm u}_n)=
T_{{\bf m},c'_1,c'_2}{\downarrow}U_q({\rm u}_n)$, then it follows
from Schur's lemma that $A$ is a constant operator on each of the subspaces
${\cal V}_{{\bf m}_n}$, on which the irreducible representations
$T_{{\bf m}_n}$ of $U_q({\rm u}_n)$ are realized. Thus, the operator $A$
is diagonal in the basis $\{ \vert {\bf m}_n,\alpha\rangle \}$
of the representation space ${\cal H}_{\bf m}$ and
$$
\langle {\bf m}_n,\alpha \vert A\vert {\bf m}_n,\alpha
\rangle =b_{{\bf m}_n},
$$
where $b_{{\bf m}_n}$ does not depend on $\alpha$. Representing the relation
(38) in the matrix form and setting $a=e_n$ and then
$a=f_n$, we obtain the set of the equalities
$$
b_{{\bf m}^{+s}_n}[l_{sn}-c_1]=b_{{\bf m}_n}[l_{sn}-c'_1] , \ \ \
b_{{\bf m}_n}[-l_{sn}+c_2]=b_{{\bf m}^{+s}_n}[-l_{sn}+c'_2] . \eqno (39)
$$
It follows from here that for $s=1,2,\cdots ,n$ and for all admissible values
of $l_{sn}$ we have
$$
\frac {[l_{sn}-c_1]}{[l_{sn}-c'_1]} =
\frac {[l_{sn}-c'_2]}{[l_{sn}-c_2]} .
$$
These equalities are possible only if
$$
c'_1=c_1+{\rm i}\pi k/h,\ \ c'_2=c_2-{\rm i}\pi k/h\ \ {\rm or\ if}\ \
c'_1=c_2+{\rm i}\pi k/h,\ \ c'_2=c_1-{\rm i}\pi k/h ,
\eqno (40)
$$
where $k\in {\bf Z}$. Thus, we proved the following
\medskip

\noindent
{\bf Proposition 2.} {\it
The irreducible representations $T_{{\bf m},c_1,c_2}$ and
$T_{{\bf m}',c'_1,c'_2}$ are equivalent if and only if ${\bf m}={\bf m}'$
and one of the conditions (40) is satisfied.}
\medskip

The recurrence equations (39) allow us to calculate the matrix
elements $b_{{\bf m}_n}\equiv b_{{\bf m}_n}({\bf m},c_1,c_2)$ of the
intertwining operator $A$ of the representations
$T_{{\bf m},c_1,c_2}$ and $T_{{\bf m},c_2,c_1}$. A direct calculation
shows that they are given by
$$
b_{{\bf m}_n}({\bf m},c_1,c_2)=b \, \frac
{a_{{\bf m}_n}({\bf m},c_1,c_2)}{a_{{\bf m}_n}({\bf m},c_2,c_1)},
$$
where $b$ is an arbitrary fixed complex number and
$a_{{\bf m}_n}({\bf m},c_1,c_2)$ is determined by one of the expressions
$$
a_{{\bf m}_n}({\bf m},c_1,c_2)=
\prod _{r=1}^k \prod _{\sigma =l_r+1}^{l_{rn}-1} [\sigma -c_2]
\prod _{s=k+1}^n \prod _{\tau =l_{sn}}^{l_{s-1}-1} [\tau -c_1] ,
\ \ \ k=1,2,\cdots ,n-1.
$$
All $n-1$ expressions for $a_{{\bf m}_n}({\bf m},c_1,c_2)$ lead (up to
a constant) to the same expression for $b_{{\bf m}_n}({\bf m},c_1,c_2)$.
Everywhere below we assume that $a_{{\bf m}_n}({\bf m},c_1,c_2)$
is determined by $k=n-1$.
\medskip

\noindent
{\bf 7. Reducible representations $T_{{\bf m},c_1,c_2}$}
\medskip

\noindent
Now we consider the representations $T_{{\bf m},c_1,c_2}$ for which the
conditions of Corollary 1 are not satisfied. Thus, in this section the
numbers $c_1$ and $c_2$ in the representations $T_{{\bf m},c_1,c_2}$
are integers.

Let $l_0,l_1,l_2,\cdots l_{n-1},l_n$ be the numbers defined by
$l_0=\infty$, $l_i=m_i-i-1$ and $l_n=-\infty$.
Let $c$ be an integer such that $l_{i-1}> c> l_i$ for some $i$
($i=1,2,\cdots ,n$).
We denote by $\Gamma ^-_c$ and $\Gamma ^+_c$
the set of the subspaces ${\cal V}_{{\bf m}_n}$
of the representation space ${\cal H}_{\bf m}$ for which the highest weights
${{\bf m}_n}$ satisfy the conditions (16) and the condition $l_{in}\equiv
m_{in}-i\le c$ and $l_{in}>c$, respectively. Let $E_c^-$ and $E^+_c$ be
the projectors mapping ${\cal H}_{\bf m}$ onto the closure of the span
of all subspaces from $\Gamma ^-_c$ and $\Gamma ^+_c$, respectively.

Under studying reducible representations $T_{{\bf m},c_1,c_2}$
we distinguish four cases depending on the placement
of $c_1$ and $c_2$ between the numbers $l_0,l_1, \cdots ,
l_{n-1},l_n$.
\medskip

{\sl Case 1:} $l_{r-1}>c_1>l_r$, $l_{s-1}>c_2>l_s$, $1\le r<s\le n$.
In this case in formula (17)
one coefficient $[l_{rn}-c_1]\omega _r({\bf m},{\bf m}_n, \alpha )$ vanishes
at $l_{rn}=c_1$ and in formula (18)
one coefficient $[-l_{sn}+c_2+1]\omega _r({\bf m},{\bf m}^{-s}_n, \alpha )$
vanishes at $l_{sn}=c_2+1$. It is seen from (17) and (18) that the operator
$T_{{\bf m},c_1,c_2}(e_n)$ only increases values of the numbers $m_{in}$
and the operator
$T_{{\bf m},c_1,c_2}(f_n)$ only decreases these values. Therefore,
vanishing of the coefficients in (17) and (18) means appearance of
invariant subspaces. There exist three invariant subspaces
$$
{\cal H}^{-+}:= E^-_{c_1}E^+_{c_2}{\cal H}_{\bf m},\ \
{\cal H}^-:= E^-_{c_1}{\cal H}_{\bf m},\ \
{\cal H}^+:=E^+_{c_2}{\cal H}_{\bf m}
$$
in the representation space ${\cal H}_{\bf m}$.
Irreducible representations of $U_q({\rm u}_{n,1})$ are realized on the
subspace ${\cal H}^{-+}$ and on the quotient spaces ${\cal H}^-/
{\cal H}^{-+}$, ${\cal H}^+/{\cal H}^{-+}$ and ${\cal H}_{\bf m}/
({\cal H}^+ +{\cal H}^-)$ (their irreducibility follows from Lemma 2).
We denote these representations by
$$
{\hat R}^{rs}({\bf m},c_1,c_2),\ \ \
R_-^{rs}({\bf m},c_1,c_2),\ \ \
R_+^{rs}({\bf m},c_1,c_2),\ \ \
{\breve R}^{rs}({\bf m},c_1,c_2),
$$
respectively.
\medskip

{\sl Case 2:} $l_{r-1}>c_1>c_2>l_r$, $1\le r\le n$. Again, vanishing of
the coefficients in (17) and (18) in this case means the representation
space ${\cal H}_{\bf m}$ has three invariant subspaces
$$
{\cal H}^{+-}:= E^-_{c_1}E^+_{c_2}{\cal H}_{\bf m},\ \ \
{\cal H}^-:= E^-_{c_1}{\cal H}_{\bf m},\ \ \
{\cal H}^+:=E^+_{c_2}{\cal H}_{\bf m}.
$$
Irreducible representations of $U_q({\rm u}_{n,1})$ are realized on the
subspace ${\cal H}^{+-}$ and on the quotient spaces ${\cal H}^-/
{\cal H}^{+-}$ and ${\cal H}^+/{\cal H}^{+-}$
(their irreducibility follows from Lemma 2).
We denote these representations respectively by
$$
{\hat R}^{rr}({\bf m},c_1,c_2),\ \ \
R_-^{rr}({\bf m},c_1,c_2),\ \ \
R_+^{rr}({\bf m},c_1,c_2).
$$

{\sl Case 3:} $c_1=l_r$, $r=1,2,\cdots ,n-1$, and $l_{s-1}>c_2>l_s$,
$1\le s\le n$. In this case the representation space has only one invariant
subspace ${\cal H}^+:=E^+_{c_2}{\cal H}_{\bf m}$. On this subspace and
on the quotient space ${\cal H}_{\bf m}/{\cal H}^+$, irreducible
representations of $U_q({\rm u}_{n,1})$ are realized. They will be denoted
respectively by
$$
{\tilde R}_+^{rs}({\bf m},c_1,c_2),\ \ \ \
{\tilde R}_-^{rs}({\bf m},c_1,c_2)
$$

{\sl Case 4:} $c_1=c_2=c$ and $l_{r-1}>c>l_r$, $1\le r\le n$. In this
case the representation $T_{{\bf m},c_1,c_2}$ is a direct sum of two
irreducible representations $R_+^{r}({\bf m},c_1,c_2)$ and
$R_-^{r}({\bf m},c_1,c_2)$ of $U_q({\rm u}_{n,1})$ which are realized
on the subspaces $E^+_{c}{\cal H}_{\bf m}$ and $E^-_{c}{\cal H}_{\bf m}$,
respectively.
\medskip

There are other two cases for consideration: the case when
$l_{r-1}>c_2>l_r$, $l_{s-1}>c_1>l_s$, $1\le r<s\le n$ and the case
when $l_{r-1}>c_2>c_1>l_r$, $1\le r\le n$. However, direct calculation
shows that the reducible representation $T_{{\bf m},c_1,c_2}$
with $l_{r-1}>c_2>l_r$, $l_{s-1}>c_1>l_s$, $1\le r<s\le n$
contains the same irreducible constituents as the representation
$T_{{\bf m},c_2,c_1}$ which belongs to Case 1 and
the reducible representation $T_{{\bf m},c_1,c_2}$
with $l_{r-1}>c_2>c_1>l_r$, $1\le r\le n$
contains the same irreducible constituents as the representation
$T_{{\bf m},c_2,c_1}$ which belongs to Case 2. The equivalence
(intertwining) operators for these irreducible constituents are
constructed in Proposition 5 below.

The results of this section show that each representation
$T_{{\bf m},c_1,c_2}$ of
$U_q({\rm u}_{n,1})$, which is not included into Corollary 1, is
reducible, that is, we have the following
\medskip

\noindent
{\bf Proposition 3.} {\it
The representation $T_{{\bf m},c_1,c_2}$
with $-\pi /2h < {\rm Im}\, c_1 \le \pi /2h$
is irreducible if and only if
$c_1$ and $c_2$ are not integers or if both $c_1$ and $c_2$
coincide with some of the numbers $l_1,l_2,\cdots , l_{n-1}$.}
\medskip

There are equivalence relations between irreducible constituents of
reducible representations of Case 1 and of Case 2. In order to find
these equivalence relations we note that equivalent irreducible
constituents must contain the same irreducible representations of the
subalgebra $U_q({\rm u}_n)$. Analysing the sets of irreducible
representations of $U_q({\rm u}_n)$
contained in irreducible constituents of reducible
representations $T_{{\bf m},c_1,c_2}$ we find pairs of
irreducible constituents which are possibly equivalent. Using the
method of section 6, we try to construct equivalence operators for these
pairs. In this way, we find that in Case 1 the irreducible representation
${\hat R}^{rs}({\bf m},c_1,c_2)$ is equivalent to
the representation $R_-^{r,s+1}({\bf m}',c'_1,c'_2)$ if $s<n$, where
the set of numbers $(l'_1,l'_2,\cdots ,l'_{n-1},c'_1,c'_2)$, $l'_j=m'_j-j-1$,
is obtained from the set of numbers $(l_1,l_2,\cdots ,l_{n-1},c_1,c_2)$,
$l_j=m_j-j-1$, by permutation of the $s$-th and $(n+1)$-th numbers,
and to the representation $R_+^{r-1,s}({\bf m}',c'_1,c'_2)$ if
$r>1$, where the set of numbers
$(l'_1,l'_2,\cdots ,l'_{n-1},c'_1,c'_2)$, $l'_j=m'_j-j-1$,
is obtained from the set of numbers $(l_1,l_2,\cdots ,l_{n-1},c_1,c_2)$,
$l_j=m_j-j-1$, by permutation of the $(r-1)$-th and $n$-th numbers.

It is similarly shown that in Case 1 the irreducible representation
${\breve R}^{rs}({\bf m},c_1,c_2)$ is equivalent to the
representation $R_-^{r+1,s}({\bf m}',c'_1,c'_2)$, where the set of numbers
$(l'_1,l'_2,\cdots$, $l'_{n-1},c'_1,c'_2)$, $l'_j=m'_j-j-1$,
is obtained from the set of numbers $(l_1,l_2,\cdots ,l_{n-1},c_1,c_2)$,
$l_j=m_j-j-1$, by permutation of the $r$-th and $n$-th numbers.

The representation
${\hat R}^{rr}({\bf m},c_1,c_2)$ from Case 2 is equivalent to the
representation $R_-^{r,r+1}({\bf m}',c'_1,c'_2)$ from Case 1 if $r<n$,
where the
set of numbers $(l'_1,l'_2,\cdots ,l'_{n-1},c'_1,c'_2)$, $l'_j=m'_j-j-1$,
is obtained from the set of numbers $(l_1,l_2,\cdots$, $l_{n-1},c_1,c_2)$,
$l_j=m_j-j-1$, by permutation of the $r$-th and $(n+1)$-th numbers
and to the representation $R_+^{r-1,r}({\bf m}',c'_1,c'_2)$ from Case 1
if $r>1$, where the
set of numbers $(l'_1,l'_2,\cdots ,l'_{n-1},c'_1,c'_2)$, $l'_j=m'_j-j-1$,
is obtained from the set of numbers $(l_1,l_2,\cdots ,l_{n-1},c_1,c_2)$,
$l_j=m_j-j-1$, by permutation of the $(r-1)$-th and $n$-th numbers.

Other equivalence relations between irreducible representations from
Cases 1--4 will be described in Theorem 2.

We also note that the irreducible representation
${\hat R}^{1,n}({\bf m},c_1,c_2)$ from Case 1 is equivalent to
the irreducible finite dimensional representation $T_{{\bf m}_{n+1}}$
of the algebra $U_q({\rm u}_{n,1})$ with highest weight
${\bf m}_{n+1}$ such that
$$
m_{1,n+1}=c_1+1,\ \ \ m_{i,n+1}=l_{i-1}+i,\  i=2,3,\cdots ,n,\ \ \
m_{n+1,n+1}=c_2+n+1.
$$

Now we can describe a structure of reducible representations
$T_{{\bf m},c_1,c_2}$. Irreducible constituents are contained in
these representations in the form of semidirect sum (that is,
there are irreducible constituents realized on invariant subspaces
and there are irreducible constituents realized on quotient spaces
which cannot be realized on invariant subspaces). Below we show
decompositions of reducible representations $T_{{\bf m},c_1,c_2}$ into
irreducible constituents. In these decompositions, two
constituents $R$ and $R'$ are connected by an arrow: $R\to R'$ if
$R'$ is a subrepresentation and $R$ is realized on a quotient space.
If constituents $R$ and $R'$ are connected as $R\oplus R'$, then they
are contained in the decomposition as a direct sum of representations,
that is, $R$ and $R'$ are realized on invariant subspaces.

For convenience, below we denote the representations $R({\bf m}, c_1,c_2)$
(with $R$ equipped with symbols) also by $R(L)$ (with $R$ equipped with
the same symbols), where $L=(l_1,l_2,\cdots ,l_{n-1}, c_1,c_2)$,
$l_i=m_i-i-1$. Then $s_{ik}L$ means $L$ with permuted $i$-th and $k$-th
numbers.

Let $c_1$ and $c_2$ be integers. If $l_{r-1}>c_1>c_2>l_r$, $1\le r\le n$,
then
$$
T_{{\bf m}, c_1,c_2}=\{ R^{rr}_-(L)\oplus R^{rr}_+(L)\} \to
R^{r,r+1}_-(s_{r,n+1}L)\ \ {\rm if} \ \ r<n, $$
$$
T_{{\bf m}, c_1,c_2}=\{ R^{rr}_-(L)\oplus R^{rr}_+(L)\} \to
R^{r-1,r}_+(s_{r-1,n}L)\ \
{\rm if} \ \ r>1.
$$
For $T_{{\bf m}, c_2,c_1}$ the decomposition coincides with the
decomposition for $T_{{\bf m}, c_1,c_2}$ if to reverse arrows to the
opposite sides.

If $l_{r-1}>c_1>l_r$, $l_{s-1}>c_2>l_s$, $1\le r<s\le n$, then
$$
T_{{\bf m}, c_1,c_2}=R_-^{r+1,s}(s_{rn}L)\to
\{ R^{rs}_-(L)\oplus R^{rs}_+(L)\} \to R^{r,s+1}_-(s_{s,n+1}L)\ \
{\rm if} \ \  s<n,  $$
$$
T_{{\bf m}, c_1,c_2}=R_-^{r+1,s}(s_{rn}L)\to
\{ R^{rs}_-(L)\oplus R^{rs}_+(L)\} \to R^{r-1,s}_+(s_{r-1,n}L)\ \
{\rm if} \ \ r>1,  $$
$$
T_{{\bf m}, c_1,c_2}=R_-^{r+1,s}(s_{rn}L)\to
\{ R^{rs}_-(L)\oplus R^{rs}_+(L)\} \to T_{{\bf m}_{n+1}}\ \
{\rm if} \ \ r=1,\ s=n,
$$
where $T_{{\bf m}_{n+1}}$ is the irreducible finite dimensional
representation of $U_q({\rm u}_{n,1})$ with highest weight ${\bf m}_{n+1}$
described above.
For $T_{{\bf m}, c_2,c_1}$ the decomposition coincides with the
decomposition for $T_{{\bf m}, c_1,c_2}$ if to reverse arrows to the
opposite sides.

If $c_1=l_r$, $r=1,2,\cdots ,n-1$, and $l_{s-1}>c_2>l_s$, $1\le s\le n$,
then
$$
T_{{\bf m}, c_1,c_2}={\tilde R}^{rs}_-(L)\to {\tilde R}^{rs}_+(L).
$$
If $c_1=c_2=c$ and $l_{r-1}>c>l_r$, $1\le r\le n$, then
$$
T_{{\bf m}, c_1,c_2}=R^r_+(L)\oplus R^r_-(L).
$$

\noindent
{\bf 8. Intertwining operators for reducible representations}
\medskip

\noindent
In section 6, we constructed intertwining operators for pairs of the
irreducible representations $T_{{\bf m},c_1,c_2}$. Intertwining operators
exist also for pair of the reducible representations. Before to construct
them we first prove the following
\medskip

\noindent
{\bf Proposition 4.} {\it Representations $T$ and $T'$ have a nonzero
intertwining operator if and only if there exist subrepresentations
$R$ and $R'$ of $T$ and $T'$, respectively, such that the quotient
representation $T/R$ and the subrepresentation $R'$ are equivalent.}
\medskip

\noindent
{\sl Proof.} Let $A$ be an intertwining operator for the representations
$T$ and $T'$, that is, $AT=T'A$. Then $T$ induces a subrepresentation $R$
on ${\rm Ker}\, A$ and $T'$ determines a subrepresentation $R'$ on
${\rm Im}\, A$. The operator $A$ induces the operator $A': {\cal H}/
{\rm Ker}\, A \to {\rm Im}\, A$ (${\cal H}$ is a representation space for
$T$) which is an equivalence operator for the representations $T/R$ and
$R'$. Conversely, let
$R$ and $R'$ be subrepresentations of $T$ and $T'$, respectively, such that
the quotient representation $T/R$ and the subrepresentation $R'$ are
equivalent and let $B$ be an equivalence operator for $T/R$ and $R'$.
We represent the representations $T$ and $T'$ in the form of semidirect
sums $T\sim {\hat R}\to R$ and $T'\sim {\hat R}'\to R'$, where
arrows mean the same as in section 7. According to these decompositions
we have the operator from the representation space of $T$ to the
representation space of $T'$ representable in the block form as
$\left( \matrix{ 0 & 0\cr
                 B & 0} \right)$.
It is an intertwining operator for the representations
$T$ and $T'$. Proposition is proved.
\medskip

In section 6, we constructed intertwining operators $A({\bf m},c_1,c_2)$
for pairs of the irreducible representations $T_{{\bf m},c_1,c_2}$
and $T_{{\bf m},c_2,c_1}$ which are diagonal in the basis
$\{ |{\bf m}_n, \alpha \rangle \}$.
Let us consider $A({\bf m},c_1,c_2)$ as an operator
function of $c_1$. Then for every fixed ${\bf m}$ and $m_0=c_1+c_2$,
$A({\bf m},c_1,c_2)$ is an analytical function of the complex variable $c_1$
in all points $c_1$ for which the representations $T_{{\bf m},c_1,c_2}$
are irreducible. This means that its diagonal matrix elements
$b_{{\bf m}_n}({\bf m},c_1,c_2)$ are analytical functions of $c_1$. Let
us take the meromorphic continuation of $A({\bf m},c_1,c_2)$ to the
points $c_1$ in which the representations $T_{{\bf m},c_1,c_2}$ are
reducible, that is, to the whole complex plane. We denote the analytically
continued operator function also by $A({\bf m},c_1,c_2)$. We
have explicit expressions for matrix elements of these operator
functions in the basis $\{ |{\bf m}_n, \alpha \rangle \}$. Analysing
analytical properties of the matrix elements
$b_{{\bf m}_n}({\bf m},c_1,c_2)$ we obtain the following
\medskip

\noindent
{\bf Proposition 5.} {\it For fixed ${\bf m}$ and $m_0=c_1+c_2$,
the operator function $A({\bf m},c_1,c_2)$ of the complex variable $c_1$
is regular in all points of the strip $-{\rm i}\pi /2h \le {\rm Im}\, c_1
\le {\rm i}\pi /2h$ except for the integral points
$c_1$ for which one of the following conditions is fulfilled:}
\medskip

(a) {\it $c_1$ and $c_2$ lie in different intervals $(l_{i-1},l_i)$ and
$(l_{j-1},l_j)$, $i=1,2,\cdots ,n-1$, $j=1,2,\cdots ,n$, and do not
coincide with their ends;}
\smallskip

(b) {\it $l_{i-1}>c_2>c_1>l_i$, $i=1,2,\cdots ,n$;}
\smallskip

(c) {\it $c_1=l_i$, $i=1,2,\cdots ,n-1$,  $l_{n-1}>c_2$;}
\smallskip

(d) {\it $c_2=l_i$, $i=1,2,\cdots ,n-1$,
$l_{j-1}>c_1>l_j$, $j=1,2,\cdots ,n-1$.}
\medskip

\noindent
{\it In the points} (a)--(d) {\it (except for the points
$l_{j-1}>c_1>l_j$, $j=1,2,\cdots ,n-1$, $c_2<l_{n-1}$ for which
$A({\bf m},c_1,c_2)$ has second order poles) the operator function
$A({\bf m},c_1,c_2)$ has first order poles.

In all regularity points except for the integral points $c_1$ for which
one of the following conditions is fulfilled:
\medskip

$(1)$ $l_{i-1}>c_1>c_2>l_i$, $i=1,2,\cdots ,n$;

$(2)$ $c_1=l_i$, $i=1,2,\cdots ,n-1$, $l_{j-1}>c_2>l_j$, $j=1,2,\cdots ,n-1;$

$(3)$ $c_2=l_i$, $i=1,2,\cdots ,n-1$, $c_1<l_{n-1}$;

$(4)$ $c_1<l_{n-1}$, $l_{i-1}>c_2>l_i$, $i=1,2,\cdots ,n-1,$
\medskip

\noindent
the kernel of the operator $A({\bf m},c_1,c_2)$ consists only of the
zero element. In the points $(1)-(4)$ this kernel coincide
\medskip

with $E^+_{c_2}E^-_{c_1}{\cal D}_{{\bf m},c_1,c_2}$ in the case $(1)$;

with $E^+_{c_2}{\cal D}_{{\bf m},c_1,c_2}$ in the case $(2)$;

with $E^-_{c_1}{\cal D}_{{\bf m},c_1,c_2}$ in the case $(3)$;

with $(1-E^+_{c_2}E^-_{c_1}){\cal D}_{{\bf m},c_1,c_2}$ in the case $(4)$,
\medskip

\noindent
where ${\cal D}_{{\bf m},c_1,c_2}$ is the span of the subspaces of
irreducible representations of $U_q({\rm u}_n)$.}

{\it In every point $c_1=c_0$, in which $A({\bf m},c_1,c_2)$ has
a first order pole, the operator
$$
B({\bf m},c_0,-c_0+m_0)={\rm Res}_{c_1=c_0} A({\bf m},c_1,c_2)
$$
is an intertwining operator for the representations
$T_{{\bf m},c_0,-c_0+m_0}$ and $T_{{\bf m},-c_0+m_0, c_0}$. If
$c_1$ satisfies the condition
$l_{i-1}>c_1\ge c_2>l_i$, $i=1,2,\cdots ,n$, then the operator
$A({\bf m},c_1,c_2)$ is a direct sum of two intertwining operators
$A({\bf m},c_1,c_2)E^+_{c_1}$ and
$A({\bf m},c_1,c_2)E^-_{c_2}$.
In every point $c_1$, in which $A({\bf m},c_1,c_2)$ has a second order pole,
the operator $B({\bf m},c_0,-c_0+m_0)$ which is a second order residue
of the operator function $A({\bf m},c_1,c_2)$ in the point $c_1=c_0$ is
an intertwining operator for the representations
$T_{{\bf m},c_0,-c_0+m_0}$ and $T_{{\bf m},-c_0+m_0, c_0}$.
The kernel of the operator $B({\bf m},c_0,-c_0+m_0)$ coincides
\medskip

with
$(1-E^+_{c_2}E^-_{c_1}){\cal D}_{{\bf m},c_1,c_2}$ if $c_1=c_0$
satisfies the condition ${\rm (a)}$;

with
$(E^-_{c_1}+E^+_{c_2}){\cal D}_{{\bf m},c_1,c_2}$ if $c_1=c_0$
satisfies the condition ${\rm (b)}$;

with
$E^+_{c_2}{\cal D}_{{\bf m},c_1,c_2}$ or
$E^-_{c_1}{\cal D}_{{\bf m},c_1,c_2}$ if $c_1=c_0$ satisfies
the condition ${\rm (c)}$ or ${\rm (d)}$, respectively.}
\medskip

Remark that Proposition 5 gives the equivalence operators for the
irreducible constituents of the reducible representations
$T_{{\bf m},c_1,c_2}$ and $T_{{\bf m},c_2,c_1}$ mentioned before
Proposition 3.
\medskip

\noindent
{\bf 9. Irreducible representations of $U_q({\rm u}_{n,1})$}
\medskip

\noindent
Now we select the set of all irreducible representations which can
be obtained from the representations $T_{{\bf m},c_1,c_2}$.
We obtain the following theorem in which ${\rm Im}\, c_1$ is considered
modulo $\pi /h$ and $c_2=m_0-c_1$.
\medskip

\noindent
{\bf Theorem 2.} {\it
The set of irreducible representations of the algebra $U_q({\rm u}_{n,1})$
consisting of all irreducible representations $T_{{\bf m},c_1,c_2}$ and all
irreducible constituents of reducible representations $T_{{\bf m},c_1,c_2}$
give the following classes of representations:}
\medskip

(a) {\it The representations $T_{{\bf m},c_1,c_2}$ for which $c_1$ and $c_2$
are not integers or for which $c_1$ and $c_2$ coincide with some of the
numbers $l_1,l_2,\cdots ,l_{n-1}$.}
\smallskip

(b) {\it The representations $R^{rs}_-({\bf m},c_1,c_2)$ and
$R^{rs}_+({\bf m},c_1,c_2)$, $1\le r\le s\le n$, $c_1>c_2$, where the numbers
$c_1$ and $c_2$ are integral no of which coincides with any of the numbers
$l_1,l_2,\cdots ,l_{n-1}$.}
\smallskip

(c) {\it The representations ${\tilde R}^{rs}_-({\bf m},c_1,c_2)$ and
${\tilde R}^{rs}_+({\bf m},c_1,c_2)$, $1\le r\le n-1$, $1\le s\le n$,
where $c_1=l_r$ and $l_{s-1}>c_2>l_s$.}
\smallskip

(d) {\it The representations $R^r_-({\bf m},c_1,c_2)$ and
$R^r_+({\bf m},c_1,c_2)$, $1\le r\le n$,
where $c_1=c_2=c$ and $l_{r-1}>c>l_r$.}
\smallskip

(e) {\it All irreducible finite dimensional representations
$T_{{\bf m}_{n+1}}$.}
\medskip

{\it Between the irreducible
representations $T_{{\bf m},c_1,c_2}$ there exist the following equivalence
relations: $T_{{\bf m},c_1,c_2}\sim T_{{\bf m},c_2,c_1}$.
In the set of the representations of classes} (b)--(e) {\it there exist
the following equivalence relations:
$$
R_-^{rs}(L)\sim R_+^{r-1,s-1}(s_{r-1,n}s_{s-1,n+1}L),\ \ \ r\ne s,\ r\ne 1,
\eqno (41)
$$
$$
R_-^{rr}(L)\sim R_+^{r-1,r-1}(s_{r-1,n}s_{n,n+1}L),\ \ \ r\ne 1,
\eqno (42)
$$
$$
{\tilde R}_-^{rs}(L)\sim {\tilde R}_+^{r,s-1}(s_{s-1,n+1}L),\ \ \ s\ne 1,
\eqno (43)
$$
$$
R_-^r(L)\sim {\tilde R}_+^{r-1,r-1}(s_{r-1,n+1}L), \ \ r\ne 1,\ \ \
R_+^r(L)\sim {\tilde R}_-^{r,r+1}(s_{r,n+1}L),\ \  r\ne n.
\eqno (44)
$$
Any another equivalence relation between the representations of classes}
(a)--(e) {\it is obtained from the above ones by means of superpositions.}
\medskip

\noindent
{\sl Proof.} The first part of the theorem follows from the above reasoning.
Existence of the equivalence relations stated in the theorem is proved by
constructing the intertwining operators for these pairs of representations
(they are diagonal in the basis $\{ \vert {\bf m}_n,\alpha\rangle \}$
and can be easily found as in section 6).

All equivalence relations in the set of the irreducible representations
$T_{{\bf m},c_1,c_2}$ were described above.
The irreducible representation $T_{{\bf m},c_1,c_2}$ cannot be equivalent
to any of the representations $R$ from classes (b)--(e) since
the restrictions $T_{{\bf m},c_1,c_2}{\downarrow} U_q({\rm u}_n)$ and
$R{\downarrow} U_q({\rm u}_n)$ do not coincide. From other side,
direct verification shows that if two representations from classes
(b)--(d) are not related by an equivalence relation from (41)--(44)
or by their superposition, then they have not coinciding restrictions
to the subalgebra $U_q({\rm u}_n)$. This proves the second part of the
theorem. Theorem is proved.
\medskip

It will be proved in a separate paper that the representations of Theorem 2
exhaust (up to tensoring by one-dimensional representations)
all irreducible representations of the algebra $U_q({\rm u}_{n,1})$.

By results of section 7, every reducible representation
$T_{{\bf m},c_1,c_2}$ has one of the forms
$$
T_{{\bf m},c_1,c_2}\sim R_1\to R_2\to R_3\to R_4 ,\ \ \
T_{{\bf m},c_1,c_2}\sim R_1\to R_2\to R_3 ,\ \ \
T_{{\bf m},c_1,c_2}\sim R_1\to R_2,
$$
where $R_i$ are irreducible constituents and in some cases an arrow
must be replaced by the symbol of a direct sum. Then the representation
$T_{{\bf m},c_2,c_1}$ is of the form
$$
T_{{\bf m},c_2,c_1}\sim R'_4\to R'_3\to R'_2\to R'_1 ,\ \ \
T_{{\bf m},c_2,c_1}\sim R'_3\to R'_2\to R'_1 ,\ \ \
T_{{\bf m},c_1,c_2}\sim R'_2\to R'_1,
$$
respectively, where $R_i$ is equivalent to $R'_i$.
Moreover, it follows from the results of section 7 that every
irreducible representation $R$ from Theorem 2 can be realized as
the constituent $R_1\sim R'_1$ in the decompositions
$$
T_{{\bf m},c_1,c_2}\sim R_1\to R_2\to \cdots \ \ \ {\rm and}\ \ \
T_{{\bf m},c_2,c_1}\sim \cdots \to R'_2\to R'_1.
$$
From other side, Proposition 5 shows that there is an intertwining operator
for $T_{{\bf m},c_1,c_2}$ and $T_{{\bf m},c_2,c_1}$ which is an equivalence
operator for the representations $R_1$ and $R'_1$, and turns into zero
on other parts of the representation space. This fact will be used in the
following section for obtaining irreducible $*$-representations of
$U_q({\rm u}_{n,1})$.
\medskip

\noindent
{\bf 10. Irreducible $*$-representations of $U_q({\rm u}_{n,1})$}
\medskip

\noindent
The aim of this section is to separate in Theorem 2 all representations
which are equivalent to $*$-representations. For this, we introduce the
notion of Hermitian-adjoint representations of $U_q({\rm u}_{n,1})$.
Two representations $T:=T_{{\bf m},c_1,c_2}$ and $T':=T_{{\bf m},c'_1,c'_2}$
on a Hilbert space ${\cal H}_{\bf m}$ are called Hermitian-adjoint
if for any $a\in U_q({\rm u}_{n,1})$ we have
$$
(T(a^*)v,v')=(v,T'(a)v')\ \ \ {\rm for\ any}\ \ \ a\in U_q({\rm u}_{n,1}).
$$
Setting here $v=v'=|{\bf m}_n,\alpha \rangle$ and $a=k_n$, we find that
$c_1+c_2=c'_1+c'_2$. Then setting
$v'=|{\bf m}_n,\alpha \rangle$, $v=|{\bf m}^{+s}_n,\alpha \rangle$
and $a=e_n$, we derive that $[l_{sn}-c_2]={\overline {[l_{sn}-c'_1]}}$,
that is, $c'_1={\overline {c_2}}+2{\rm i}\pi k/h$ and
$c'_2={\overline {c_1}}-2{\rm i}\pi k/h$, where $k\in {\bf Z}$.
Direct calculation shows that the
representations $T:=T_{{\bf m},c_1,c_2}$ and
$T':=T_{{\bf m}, {\overline {c_2}},{\overline {c_1}}}$
are indeed Hermitian-adjoint.
Thus, if $-\pi /2h < {\rm Im}\, c_1,{\rm Im}\, c'_1\le \pi /2h$, then
{\it the representations
$T:=T_{{\bf m},c_1,c_2}$ and $T':=T_{{\bf m},c'_1,c'_2}$ are
Hermitian-adjoint if and only if
$c'_1={\overline {c_2}}$ and $c'_2={\overline {c_1}}$.}

Now we shall find which of the representations of Theorem 2 are equivalent
to $*$-representations. The representation $T:=T_{{\bf m},c_1,c_2}$
on the Hilbert space ${\cal H}_{\bf m}$ is equivalent to such a
representation if there exists a scalar product (that is, a strictly
positive Hermitian form) $H(\cdot ,\cdot )$ of ${\cal D}_{\bf m}$ such
that
$$
H(T(a^*)v,v')=H(v,T(a)v'),\ \ \  a\in U_q({\rm u}_{n,1}). \eqno (45)
$$
Any Hermitian form $H(\cdot ,\cdot )$ on ${\cal D}_{\bf m}$ can be
represented as
$$
H(v,v')=(v,Qv'),
$$
where $(\cdot ,\cdot )$ is the scalar product on ${\cal H}_{\bf m}$
and $Q$ is a Hermitian operator on ${\cal D}_{\bf m}$.
It follows from (45) that
$$
H(v,T(a)v') =(v,QT(a)v')=H(T(a^*)v,v')=(T(a^*)v,Qv')=(v,T'(a)Qv'),
$$
where $T':=T_{{\bf m},{\overline {c_2}},{\overline {c_1}}}$. Hence,
$QT(a)=T'(a)Q$, that is, $Q$ is an intertwining operator for the
representations $T=T_{{\bf m},c_1,c_2}$ and
$T'=T_{{\bf m},{\overline {c_2}},{\overline {c_1}}}$. As we have seen in
section 6, if the representation
$T$ is irreducible, then such an operator exists if
${\overline {c_2}}=c_1$ or
${\overline {c_2}}=c_2$, ${\overline {c_1}}=c_1$ or
${\overline {c_2}}=c_2+{\rm i}\pi /h$, ${\overline {c_1}}=c_1-{\rm i}\pi /h$.
In the first case the representation $T=T_{{\bf m},c_1,c_2}$ coincides
with $T'=T_{{\bf m},{\overline {c_2}},{\overline {c_1}}}$ and
the operator $Q$ is multiple to the identity operator. In the second
case, the numbers $c_1$ and $c_2$ are real, and the operator
$Q$ is the intertwining operator for $T=T_{{\bf m},c_1,c_2}$
and $T'=T_{{\bf m},c_2,c_1}$ constructed in section 6.
In the third case ${\rm Im}\, c_1=-{\rm Im}\, c_2=\pi /2h$ and the
intertwining operator is from section 6.

In the first case, the Hermitian form $H(\cdot ,\cdot )$ coincides with
the scalar product and, consequently, the representation
$T_{{\bf m},c_1,c_2}$ with ${\overline {c_2}}=c_1$ is a $*$-representation
with respect to the scalar product $(\cdot ,\cdot)$.
In the second and third cases, the Hermitian form has the form
$H(v,v')=(v, Qv')$,
where $Q$ is a non-trivial intertwining operator. We have to find
when this form is positive.
Since the matrices of all constructed intertwining operators are
diagonal with respect to the basis $\{ | {\bf m}_n,\alpha\rangle \}$,
then we must find when matrix elements of the intertwining operator
are all positive or all negative.
Making these calculations explicitly we obtain all irreducible
representations $T_{{\bf m},c_1,c_2}$ which are equivalent to
$*$-representations.

The same method is appropriate for the irreducible representations
from classes (b)--(e) of Theorem 2. For every representation $R$
from this set there exist (according to the remark at the end of section 9)
the representations $T_{{\bf m},c_1,c_2}$ and $T_{{\bf m},c_2,c_1}$
such that $T_{{\bf m},c_1,c_2}\sim R\to \cdots $ and
$T_{{\bf m},c_2,c_1}\sim \cdots \to R'$, $R\sim R'$.
Now we continue as above. The difference is that now
the Hermitian form $H(v,v')=(v,Qv')$ must be only positive (not
strictly positive). As above, we derive that $QT_{{\bf m},c_1,c_2}(a)=
T_{{\bf m},c_2,c_1}(a)Q$, $a\in U_q({\rm u}_{n,1})$, and $Q$ is an
equivalence operator for $R$ and $R'$. On the carrier space ${\cal H}$
of $R'$ the Hermitian form $H$ is strictly positive if $R$ is
equivalent to a $*$-representation. The form $H$ is strictly positive
on ${\cal H}$ if $Q$ is a positive (or negative) Hermitian operator, that
is, if his non-vanishing diagonal matrix elements are all positive or all
negative. Making these calculations explicitly, we
obtain the theorem formulated below. In the formulation of this theorem,
a series of numbers $a_1,a_2,\cdots ,a_k$
is called {\it dense} if $a_i=a_{i-1}-1$, $i=2,3,\cdots ,k$.
\medskip

\noindent
{\bf Theorem 3.} {\it
The following irreducible representations of $U_q({\rm u}_{n,1})$ from
Theorem 2 are equivalent to $*$-representations:}
\medskip

(a) {\it the representations $T_{{\bf m},c_1,c_2}$, $c_1=\overline {c_2}$
(principal series of $*$-representations);}
\smallskip

(b) {\it the representations $T_{{\bf m},c_1,c_2}$,
${\rm Im}\, c_1=-{\rm Im}\, c_2=\pi /2h$ (the strange series);}
\smallskip

(c) {\it the representations $T_{{\bf m},c_1,c_2}$, where $c_1$ and $c_2$ are
real numbers for which there exist numbers $l_r=m_r-r-1$, $l_s=m_s-s-1$,
$r,s=1,2,\cdots ,n-1$, such that $\vert l_r-c_1\vert <1$, $\vert l_s-c_2\vert
<1$ and the series $l_r,l_{r+1},\cdots l_s$ (for $c_1>c_2$) or the series
$l_s,l_{s+1},\cdots , l_r$ (for $c_1<c_2$) is dense (the supplementary
series);}
\smallskip

(d) {\it the representations $R^{ij}_-({\bf m},c_1,c_2)$, $c_1>c_2$, if
$i=j$ or if the series $c_1,l_i,l_{i+1},\cdots ,$ $l_{j-1}$ is dense;}
\smallskip

(e) {\it the representations $R^{ij}_+({\bf m},c_1,c_2)$, $c_1>c_2$, if
$i=j$ or if the series $l_i,l_{i+1},\cdots ,$ $l_{j-1},c_2$ is dense;}
\smallskip

(f) {\it the representations ${\tilde R}^{ij}_+({\bf m},c_1,c_2)$, if
$i<j$ and the series $l_i,l_{i+1},\cdots ,l_{j-1},c_2$ is dense or if
$i\ge j$ and the series $l_j,l_{j+1},\cdots ,l_i$ is dense;}
\smallskip

(g) {\it the representations ${\tilde R}^{ij}_-({\bf m},c_1,c_2)$, if
$i<j$ and the series $l_i,l_{i+1},\cdots ,l_{j-1}$ is dense or if
$i\ge j$ and the series $c_2,l_j,l_{j+1},\cdots ,l_i$ is dense;}
\smallskip

(h) {\it all the representations $R^{i}_+({\bf m},c,c)$ and
$R^{i}_-({\bf m},c,c)$.}
\medskip

Note that there is a one-to-one
correspondence between nonequivalent irreducible $*$-representations
of $U_q({\rm u}_{n,1})$ which are irreducible components of reducible
representations $T({\bf m},c_1,c_2)$ and nonequivalent unitary
irreducible representations of the same type of the group $U(n,1)$.
The list of the last representations can be found in [14].

Remark that the strange series of representations disappears when $q\to 1$.
This means that this series is absent for the classical case.

Theorem 3 classifies irreducible $*$-representations of the algebra
$U_q({\rm u}_{n,1})$ in the set of representations of Theorem 2.
However, the algebra $U_q({\rm u}_{n,1})$ has also
irreducible $*$-representations belonging to the set of
representations obtained from the representations of Theorem 2 by
tensoring by one-dimensional
representations. Below, we give a list of such representations.
\medskip

\noindent
{\bf Proposition 6.} (a) {\it If $T$ is a $*$-representation from Theorem 3,
then the representation $T'$ determined by the operators $T'(k_j)=-T(k_j)$,
$T'(e_r)=T(e_r)$, $T'(f_r)=T(f_r)$, $j=1,2,\cdots ,n+1,$ $r=1,2,\cdots ,n,$
is a $*$-representation of $U_q({\rm u}_{n,1})$.}
\smallskip

(b) {\it Let $T_{{\bf m}_{n+1}}$ be a finite dimensional representation
of $U_q({\rm u}_{n,1})$ defined by formulas (12)--(14), written down
for $U_q({\rm gl}_{n+1})$. Then the representation $T'_{{\bf m}_{n+1}}$,
determined by the formulas
$T'_{{\bf m}_{n+1}}(k_{n+1})=-T_{{\bf m}_{n+1}}(k_{n+1})$,
$T'_{{\bf m}_{n+1}}(e_n)=T_{{\bf m}_{n+1}}(e_n)$,
$T'_{{\bf m}_{n+1}}(f_n)=-T_{{\bf m}_{n+1}}(f_n)$ and
$T'_{{\bf m}_{n+1}}(a)=T_{{\bf m}_{n+1}}(a)$ for all other
generating elements of $U_q({\rm gl}_{n+1})$, is a $*$-representation of
$U_q({\rm u}_{n,1})$.}
\medskip

\noindent
{\sl Proof} is given by a direct verification..
\medskip

By Proposition 6, $U_q({\rm u}_{n,1})$ has finite
dimensional irreducible $*$-representations. It is not a case for
the classical case.
\bigskip

The research described in this publication was made possible in part by
Award No. UP1-309 of the U.S. Civilian Research and Development
Foundation for the Independent States of the Former Soviet Union (CRDF),
and by Ukrainian DFFD Grant 1.4/206.

\end{document}